\documentclass[12pt]{article}
\usepackage{latexsym, amssymb}
\DeclareMathAlphabet\gothic{U}{euf}{m}{n}

\makeatletter
\def\eqnarray{\stepcounter{equation}\let\@currentlabel=\theequation
\global\@eqnswtrue
\tabskip\@centering\let\\=\@eqncr
$$\halign to \displaywidth\bgroup\hfil\global\@eqcnt\z@
  $\displaystyle\tabskip\z@{##}$&\global\@eqcnt\@ne
  \hfil$\displaystyle{{}##{}}$\hfil
  &\global\@eqcnt\tw@ $\displaystyle{##}$\hfil
  \tabskip\@centering&\llap{##}\tabskip\z@\cr}

\def\endeqnarray{\@@eqncr\egroup
      \global\advance\c@equation\m@ne$$\global\@ignoretrue}
\makeatother

\begin{document}
\bibliographystyle{tom}

\newtheorem{lemma}{Lemma}[section]
\newtheorem{thm}[lemma]{Theorem}
\newtheorem{cor}[lemma]{Corollary}
\newtheorem{voorb}[lemma]{Example}
\newtheorem{rem}[lemma]{Remark}
\newtheorem{prop}[lemma]{Proposition}
\newtheorem{stat}[lemma]{{\hspace{-5pt}}}
\newtheorem{obs}[lemma]{Observation}

\newenvironment{remarkn}{\begin{rem} \rm}{\end{rem}}
\newenvironment{exam}{\begin{voorb} \rm}{\end{voorb}}
\newtheorem{rems}[lemma]{Remarks}
\newenvironment{remarkns}{\begin{rems} \rm}{\end{rems}}

\newenvironment{obsn}{\begin{obs} \rm}{\end{obs}}
\newenvironment{obsns}{\begin{obss} \rm}{\end{obss}}

\newtheorem{obss}[lemma]{Observations}

\newcommand{\gota}{\gothic{a}}
\newcommand{\gotb}{\gothic{b}}
\newcommand{\gotc}{\gothic{c}}
\newcommand{\gote}{\gothic{e}}
\newcommand{\gotf}{\gothic{f}}
\newcommand{\gotg}{\gothic{g}}
\newcommand{\gothh}{\gothic{h}}
\newcommand{\gotk}{\gothic{k}}
\newcommand{\gotm}{\gothic{m}}
\newcommand{\gotn}{\gothic{n}}
\newcommand{\gotp}{\gothic{p}}
\newcommand{\gotq}{\gothic{q}}
\newcommand{\gotr}{\gothic{r}}
\newcommand{\gots}{\gothic{s}}
\newcommand{\gotu}{\gothic{u}}
\newcommand{\gotv}{\gothic{v}}
\newcommand{\gotw}{\gothic{w}}
\newcommand{\gotz}{\gothic{z}}
\newcommand{\gotG}{\gothic{G}}
\newcommand{\gotL}{\gothic{L}}
\newcommand{\gotS}{\gothic{S}}
\newcommand{\gotT}{\gothic{T}}

\newcounter{teller}
\renewcommand{\theteller}{\Roman{teller}}
\newenvironment{tabel}{\begin{list}%
{\rm \bf \Roman{teller}.\hfill}{\usecounter{teller} \leftmargin=1.1cm
\labelwidth=1.1cm \labelsep=0cm \parsep=0cm}
                      }{\end{list}}

\newcounter{tellerr}
\renewcommand{\thetellerr}{\roman{tellerr}}
\newenvironment{subtabel}{\begin{list}%
{\rm  \roman{tellerr}.\hfill}{\usecounter{tellerr} \leftmargin=1.1cm
\labelwidth=1.1cm \labelsep=0cm \parsep=0cm}
                         }{\end{list}}

\newcounter{proofstep}
\newcommand{\nextstep}{\refstepcounter{proofstep}\ruimte \par 
          \noindent{\bf Step \theproofstep} \hspace{5pt}}
\newcommand{\firststep}{\setcounter{proofstep}{0}\nextstep}

\newcommand{\Ni}{{\bf N}}
\newcommand{\Ri}{{\bf R}}
\newcommand{\Ci}{{\bf C}}
\newcommand{\Ti}{{\bf T}}
\newcommand{\Zi}{{\bf Z}}
\newcommand{\Fi}{{\bf F}}

\newcommand{\proof}{\mbox{\bf Proof} \hspace{5pt}} 
\newcommand{\remark}{\mbox{\bf Remark} \hspace{5pt}}
\newcommand{\ruimte}{\vskip10.0pt plus 4.0pt minus 6.0pt}

\newcommand{\ad}{{\mathop{\rm ad}}}
\newcommand{\Ad}{{\mathop{\rm Ad}}}
\newcommand{\Aut}{\mathop{\rm Aut}}
\newcommand{\arccot}{\mathop{\rm arccot}}
\newcommand{\clos}{{\mathop{\rm cl}}}

\newcommand{\diam}{\mathop{\rm diam}}
\newcommand{\divv}{\mathop{\rm div}}
\newcommand{\codim}{\mathop{\rm codim}}
\newcommand{\interior}{\mathop{\rm int}}
\newcommand{\RRe}{\mathop{\rm Re}}
\newcommand{\IIm}{\mathop{\rm Im}}
\newcommand{\Tr}{{\mathop{\rm Tr}}}
\newcommand{\supp}{\mathop{\rm supp}}
\newcommand{\sgn}{\mathop{\rm sgn}}
\newcommand{\esssup}{\mathop{\rm ess\,sup}}
\newcommand{\Int}{\mathop{\rm Int}}
\newcommand{\Leibniz}{\mathop{\rm Leibniz}}
\newcommand{\lcm}{\mathop{\rm lcm}}
\newcommand{\mod}{\mathop{\rm mod}}
\newcommand{\spann}{\mathop{\rm span}}
\newcommand{\ubar}{\underline{\;}}
\newcommand{\one}{1\hspace{-4.5pt}1}

\hyphenation{groups}
\hyphenation{unitary}

\newcommand{\dm}{{d_{\mbox{\gothicss m}}}}
\newcommand{\dms}{{d_{\mbox{\gothics m}}}}
\newcommand{\hc}{\overline{H}}
\newcommand{\sodp}{ \{ 1,\ldots,d' \} }
\newcommand{\ccig}{C_c^\infty(G)}
\newcommand{\Cprime}[1]{C^{#1}{\hspace{0.8pt}}^\prime}
\newcommand{\CCprime}[1]{{\cal C}^{#1}{\hspace{0.8pt}}^\prime}
\newcommand{\Cdeltaprime}[1]{C_\Delta^{#1}{\hspace{0.5pt}}^\prime}

\newcommand{\Gdual}{\widetilde{G}}
\newcommand{\htozp}{H'_{2;1}\raisebox{10pt}[0pt][0pt]{\makebox[0pt]{\hspace{-28pt}$\scriptstyle\circ$}}}
\newcommand{\htozpu}{H^{\prime(u)}_{2;1}\raisebox{10pt}[0pt][0pt]{\makebox[0pt]{\hspace{-36pt}$\scriptstyle\circ$}}}

\newcommand{\htozpud}{H^{(u)}_{2;1}\raisebox{10pt}[0pt][0pt]{\makebox[0pt]{\hspace{-32pt}$\scriptstyle\circ$}}}
\newcommand{\htozpukd}{H^{(u_k)}_{2;1}\raisebox{10pt}[0pt][0pt]{\makebox[0pt]{\hspace{-41pt}$\scriptstyle\circ$}}}
\newcommand{\htozpund}{H^{(u_n)}_{2;1}\raisebox{10pt}[0pt][0pt]{\makebox[0pt]{\hspace{-41pt}$\scriptstyle\circ$}}}

\newcommand{\htozpH}{H^{\prime(H)}_{2;1}\raisebox{10pt}[0pt][0pt]{\makebox[0pt]{\hspace{-36pt}$\scriptstyle\circ$}}}
\newcommand{\htozpun}{H^{\prime(u_n)}_{2;1}\raisebox{10pt}[0pt][0pt]{\makebox[0pt]{\hspace{-45pt}$\scriptstyle\circ$}}}
\newcommand{\htoz}{H_{2;1}\raisebox{10pt}[0pt][0pt]{\makebox[0pt]{\hspace{-28pt}$\scriptstyle\circ$}}}
\newcommand{\twocases}{Suppose
   \begin{list}{}{\leftmargin=40mm \labelwidth=20mm \labelsep=5mm}
    \item[{\hfill\rm \bf either}\hfill] $H$ is  strongly elliptic 
    \item[{\hfill\rm \bf or}\hfill] $G$ is stratified and 
          $a_1,\ldots,a_{d'}$ is a basis for $\gotg_1$ 
     \newline
   in the stratification $(\gotg_m)_{m \in \{ 1,\ldots,r \} }$ of $\gotg$.
        \end{list}    }

\newcommand{\lstar}[1]{\mathbin{{}_{#1}\mskip-1mu*}}
\newcommand{\ldoublestar}[2]{\mathbin{{}_{#2}\mskip-1mu({}_{#1}\mskip-1mu*)}}
\newcommand{\sdp}[1]{\mathbin{{}_{#1}\mskip-2mu\times}}
\newcommand{\tfrac}[2]{{\textstyle \frac{#1}{#2} }}
\newcommand{\sptilde}{\hspace{0pt} \widetilde{\rule{0pt}{8pt} \hspace{7pt} }}

\newcommand{\cc}{{\cal C}}
\newcommand{\cd}{{\cal D}}
\newcommand{\cf}{{\cal F}}
\newcommand{\ch}{{\cal H}}
\newcommand{\ck}{{\cal K}}
\newcommand{\cl}{{\cal L}}
\newcommand{\cm}{{\cal M}}
\newcommand{\co}{{\cal O}}
\newcommand{\cs}{{\cal S}}
\newcommand{\ct}{{\cal T}}
\newcommand{\cx}{{\cal X}}
\newcommand{\cy}{{\cal Y}}
\newcommand{\cz}{{\cal Z}}

\newfont{\fontcmrten}{cmr10}
\newcommand{\slbrl}{\mbox{\fontcmrten (}}
\newcommand{\slbrr}{\mbox{\fontcmrten )}}

\newlength{\hightcharacter}
\newlength{\widthcharacter}
\newcommand{\covsup}[1]{\settowidth{\widthcharacter}{$#1$}\addtolength{\widthcharacter}{-0.15em}\settoheight{\hightcharacter}{$#1$}\addtolength{\hightcharacter}{0.1ex}#1\raisebox{\hightcharacter}[0pt][0pt]{\makebox[0pt]{\hspace{-\widthcharacter}$\scriptstyle\circ$}}}
\newcommand{\cov}[1]{\settowidth{\widthcharacter}{$#1$}\addtolength{\widthcharacter}{-0.15em}\settoheight{\hightcharacter}{$#1$}\addtolength{\hightcharacter}{0.1ex}#1\raisebox{\hightcharacter}{\makebox[0pt]{\hspace{-\widthcharacter}$\scriptstyle\circ$}}}
\newcommand{\scov}[1]{\settowidth{\widthcharacter}{$#1$}\addtolength{\widthcharacter}{-0.15em}\settoheight{\hightcharacter}{$#1$}\addtolength{\hightcharacter}{0.1ex}#1\raisebox{0.7\hightcharacter}{\makebox[0pt]{\hspace{-\widthcharacter}$\scriptstyle\circ$}}}
\newcommand{\gznp}{\settowidth{\widthcharacter}{$G$}\addtolength{\widthcharacter}{-0.15em}\settoheight{\hightcharacter}{$G$}\addtolength{\hightcharacter}{0.1ex}G\raisebox{\hightcharacter}[0pt][0pt]{\makebox[0pt]{\hspace{-\widthcharacter}$\scriptstyle\circ$}}_N'}
\newcommand{\sgznp}{{\settowidth{\widthcharacter}{$G$}\addtolength{\widthcharacter}{-0.15em}\settoheight{\hightcharacter}{$G$}\addtolength{\hightcharacter}{0.1ex}\raisebox{-5pt}{$\scriptstyle G\raisebox{0.7\hightcharacter}[0pt][0pt]{\makebox[0pt]{\hspace{-\widthcharacter}$\scriptstyle\circ$}}_N'$}}}

 \thispagestyle{empty}

\begin{center}
{\Large{\bf Analysis of degenerate elliptic operators  }}\\[2mm] 
{\Large{\bf of  Gru\v{s}in type}}  \\[2mm]
\large  Derek W. Robinson$^1$ and Adam Sikora$^2$\\[2mm]

\normalsize{January 2013}
\end{center}

\vspace{5mm}

\begin{center}
{\bf Abstract}
\end{center}

\begin{list}{}{\leftmargin=1.8cm \rightmargin=1.8cm \listparindent=10mm 
   \parsep=0pt}
   \item
We analyze degenerate, second-order, elliptic operators $H$ in divergence form on $L_2(\Ri^{n}\times\Ri^{m})$.
We assume the coefficients are real symmetric and $a_1H_\delta\geq H\geq a_2H_\delta$ for some
$a_1,a_2>0$ where
\[
H_\delta=-{\nabla}_{x_1}\cdot(c_{\delta_1, \delta'_1}(x_1)\,\nabla_{x_1})-c_{\delta_2, \delta'_2}(x_1)\,\nabla_{x_2}^2
\;.
\]
Here $x_1\in\Ri^n$, $x_2\in\Ri^m$ and $c_{\delta_i, \delta'_i}$ are positive measurable functions
such that $c_{\delta_i, \delta'_i}(x)$ behaves like $|x|^{\delta_i}$ as $x\to0$ and 
$|x|^{\delta_i'}$ as $x\to\infty$ with $\delta_1,\delta_1'\in[0,1\rangle$
 and  $\delta_2,\delta_2'\geq0$.

Our principal results state that the submarkovian semigroup $S_t=e^{-tH}$ is conservative and its
kernel $K_t$ satisfies  bounds
\[
0\leq  K_t(x\,;y)\leq a\,(|B(x\,;t^{1/2})|\,|B(y\,;t^{1/2})|)^{-1/2}
\]
where  $|B(x\,;r)|$ denotes the volume of the ball $B(x\,;r)$  centred at $x$ with radius $r$
measured with respect to the  Riemannian distance associated with $H$.
The proofs depend on detailed subelliptic estimations on $H$, a precise characterization of the 
Riemannian distance and the corresponding volumes and wave equation techniques which
exploit the finite speed of propagation.

We discuss further implications of these bounds and give explicit examples that show the kernel is not necessarily strictly positive, nor continuous.

\end{list}

\vspace{1cm}


\vspace{1cm}
%


\vspace{1cm}

\noindent
\begin{tabular}{@{}cl@{\hspace{10mm}}cl}
1. &Mathematical Sciences Institute  & 
  2. &Department of Mathematical Sciences  \\
& Australian National University  & 
  &   New Mexico State University  \\
& Canberra, ACT 0200& 
  &  P.O. Box 30001 \\
& Australia.& 
  & Las Cruces  \\
&{}& 
  & NM 88003-8001, USA.
\end{tabular}

\newpage
\setcounter{page}{1}

\section{Introduction}\label{Scsg1}

Our aim is to analyze solutions of the parabolic evolution equations associated with degenerate
elliptic,   second-order  operators,   $H$ in divergence form
on $L_2(\Ri^d)$, i.e.\  operators formally expressed as 
\begin{equation}
H=-\sum^d_{i,j=1}\partial_i\,c_{ij}\,\partial_j
\;,
\label{ecsg1.1}
\end{equation}
where $\partial_i=\partial/\partial x_i$, the $c_{ij}$ are real-valued measurable functions
and  the  coefficient matrix $C=(c_{ij})$ is symmetric and  positive-definite almost-everywhere.
We assume that $d=n+m$ and $C\sim C_\delta$ where $C_\delta$  is  a block diagonal matrix, 
$C_\delta(x_1,x_2)=c_{\delta_1, \delta'_1}(x_1)\, I_n+c_{\delta_2, \delta'_2}(x_1)\,I_m$,
on $\Ri^n\times\Ri^m$
with  $c_{\delta_1, \delta_1'},c_{\delta_2, \delta'_2}$   positive measurable  functions and
\begin{equation}
c_{\delta_i, \delta_i'}(x)\sim |x|^{(2\delta_i,2\delta_i')}
\;.
\label{ecsg1.10}
\end{equation}
The indices $\delta_1,\delta_2,\delta'_1,\delta'_2$ are all non-negative and $\delta_1,\delta_1'<1$
but there is no upper bound on $\delta_2$ and $\delta'_2$.
(Here and in the sequel we use the notation  $f\sim g$ for two functions $f$, $g$ with values in an ordered space   if there are $a$, $a'>0$ such that $a\,f\leq g\leq a'\,f$ uniformly.
In addition $a^{(\alpha,\alpha')}= a^\alpha$ if $ a \leq 1 $
and  $a^{(\alpha,\alpha')}=  a^{\alpha'}$ if $   a \geq 1  $.)
Then  $H\sim H_\delta$
where
\begin{equation}
H_\delta=-{\nabla}_{x_1}\cdot(c_{\delta_1, \delta'_1}\,\nabla_{x_1})-c_{\delta_2, \delta'_2}\,\nabla_{x_2}^2
\label{ecsg1.2}
\end{equation}
with  $x_1\in\Ri^n$, $x_2\in\Ri^m$, $\nabla_{x_1}$,   $\nabla_{x_2}$ the gradient operators on $L_2(\Ri^n)$ and   $L_2(\Ri^m)$, respectively.
The operators $H$ and $H_\delta$ will be defined precisely in terms of Dirichlet forms $h$ and $h_\delta$ in Section~\ref{Scsg2} with  $h\sim h_\delta$ in the sense of ordering of quadratic forms.
Note  we do not assume any regularity of the coefficients $C=(c_{ij})$ but since $C$ is only defined up to equivalence with $C_\delta$ there is a freedom of choice of the latter coefficients.
In particular they may be chosen to be continuous.

\smallskip

We refer to the foregoing  operators as Gru\v{s}in operators since the  
 $H_\delta$  are superficially similar  to the 
two-dimensional degenerate operators $H_k= -\partial_1^2-x_1^{2k}\,\partial_2^2$, with $k\in\Ni$, 
  introduced in \cite{Gru}.
The properties of the $H_\delta$ with $n=1$ can, however,  differ dramatically
  from those of the $H_k$.
  The degeneracy of the $H_k$ on the hyperplane $x_1=0$ is due to a degeneracy of the tangential 
  component $x_1^k\,\partial_2$ of the underlying flow.
  But the $H_\delta$, with $n=1$, can have a degeneracy in the normal component of the flow.
If this latter degeneracy is sufficiently strong, e.g.\ if $\delta_1\in[1/2,1\rangle$,  then the 
the corresponding heat kernel is neither strictly positive nor continuous since the
 evolution leaves the two  half-spaces $x_1\geq0$ and $x_1\leq 0$ invariant.
This non-ergodic behaviour  is a reflection of the properties of the one-dimensional operator $H=-\partial\,|x|\,\partial$ (see \cite{ERSZ1}, Sections~5 and 6).
In this example the  Riemannian distance, $d(x\,;y)=2\,|x-y|^{1/2}$, is well-defined but the corresponding diffusion is not ergodic
This illustrates that  the Riemannian geometry is not necessarily appropriate for the 
description of properties of  strongly  degenerate operators.
Indeed the small time asymptotics of the kernel is usually given by a larger distance \cite{HiR}
\cite{AH}   \cite{ERS4} which incorporates possible separation phenomena.

\smallskip

Since the theory of elliptic operators is such a vast subject it is not possible to cite all relevant material.
Background information on the theory of degenerate elliptic operators can be found in the books and reviews \cite{GT}  \cite{JSC1} \cite{Haj} \cite{Sal5}  \cite{SawW} and references therein.
More specifically \cite{Tru2} \cite{FKS} \cite{FL} \cite{FL2} \cite{FrS} \cite{FGW2} analyze various classes of degenerate elliptic operators with $L_\infty$-coefficients principally with a view to establishing the main regularity properties 
of the strongly elliptic theory, e.g.\ H\"older regularity of solutions, Harnack inequalities, etc.
In particular \cite{FL2} \cite{FrS} and  \cite{FGW2} deal with generalizations of the Gru\v{s}in operator (\ref{ecsg1.2})
under different assumptions on the coefficients.
The situation of $n=1$ and $\delta_1\in[1/2,1\rangle$, i.e.\ the situation for which the solutions are not continuous and the 
evolution not ergodic, is not covered by these authors.


\section{Preliminaries}\label{Scsg2}

We begin with the precise definition of  the degenerate  elliptic operators formally given by   (\ref{ecsg1.1}) or  (\ref{ecsg1.2}).
If $C=(c_{ij})$ is  a real-valued symmetric matrix with measurable, locally integrable, coefficients $c_{ij}$
which is positive-definite almost everywhere one can
define  the   positive quadratic  form $h$  by 
\begin{equation}
h(\varphi)=\sum^d_{i,j=1}(\partial_i\varphi,c_{ij}\partial_j\varphi)
= \int_{\Ri^d}dx\,\sum^d_{i,j=1}c_{ij}(x)(\partial_i\varphi)(x)(\partial_j\varphi)(x)
\label{ecsg2.1}
\end{equation}
for $\varphi\in D(h)$ where $D(h)$ consists 
of those $\varphi\in W^{1,2}(\Ri^d)$ for which the integral converges.
The form is   not  necessarily closed nor even  closable.
The key property is the following.

\begin{lemma}\label{lcsg1.2}
The form $h$ of a Gru\v{s}in operator is closable.
The closure is a Dirichlet form.
\end{lemma}
\proof\
It follows by definition that  $h\sim h_\delta$, in the sense of quadratic forms, where
\begin{equation}
h_\delta(\varphi)=(\nabla_{x_1}\varphi,c_{\delta_1, \delta'_1}\,\nabla_{x_1}\varphi)
+(\nabla_{x_2}\varphi,c_{\delta_2, \delta'_2}\,\nabla_{x_2}\varphi)
\label{ecsg2.10}
\end{equation}
for all $\varphi\in D(h_\delta)=D(h)$.
Therefore  it suffices to prove that $h_\delta$ is closable.
But we may assume the coefficients  $c_{\delta_1, \delta_1'},c_{\delta_2, \delta'_2}$ are continuous.
Then the closability of $h_\delta$ follows by standard reasoning (see, for example,  \cite{MR}, Section~II.2b).
Finally the closure $\overline h$ of $h$ is a Dirichlet form by a `viscosity' approximation argument
(see Section~2 of \cite{ERSZ1}).
\hfill$\Box$

\bigskip

We now define the Gru\v{s}in operator to be the positive self-adjoint operator $H$ associated with the closure $\overline h$ of the quadratic form $h$.
Further $S$ will denote the submarkovian semigroup generated by $H$ and $K$ the corresponding non-negative integral kernel.

\smallskip

Our   aim in Section~\ref{Scsg3}  is to obtain {\it a priori} bounds on $t\to K_t$  uniform over $\Ri^d$.
 This is equivalent to obtaining bounds  on the crossnorms 
 $\|S_t\|_{2\to\infty}$ of the semigroup 
 as a map from $L_2$ to $L_\infty$  since 
   \[
\esssup_{x, y\in\Ri^d}|K_t(x\,;y)|
 =\esssup_{x\in\Ri^d}\int_{\Ri^d}dy\,|K_{t/2}(x\,;y)|^2=(\|S_{t/2}\|_{2\to\infty})^2
 \]
 for all $t>0$.
 Then $S_t$ is automatically bounded as an operator from $L_1$ to $L_\infty$
 since $\|S_t\|_{1\to\infty}=(\|S_{t/2}\|_{2\to\infty})^2$.
 The standard method of obtaining bounds on the crossnorms $\|S_t\|_{2\to\infty}$ is via Nash inequalities and for the Gru\v{s}in operators it is necessary to  consider  inequalities  involving Fourier  multipliers.

 \smallskip
 
 Let $ F$ be a positive real function over $\Ri^d$ and
 define the positive self-adjoint operator $F=F(i\nabla_x)$  on 
$L_2(\Ri^d)$ by spectral theory.
 The operator $F$ acts by multiplication of the Fourier transform.
 Specifically $\widetilde{(F\varphi)}(p)= F(p)\widetilde\varphi(p)$
where  $\widetilde\varphi$ denotes
the Fourier transform of $\varphi\in L_2(\Ri^d)$.
Next  let $f $ denote the closed form  corresponding to $F $, i.e.\
\[
f(\varphi)=\int_{\Ri^d}dp\,F(p)\,|\widetilde\varphi(p)|^2
\]
with $D(f)$ the subspace of $\varphi\in L_2(\Ri^d)$ for which the integral is finite.
Finally  let $V_F(r)$ denote the volume (Lebesgue measure) of the set $\{p: F(p)<r^2\}$.

In the next three lemmas $h$ is a general Dirichlet  form on $L_2(\Ri^d)$.
 The estimates are not restricted to Gru\v{s}in operators.

\begin{lemma} \label{lcsg2.1}
If  $h\geq f$ then
\begin{equation}
\|\varphi\|_2^2\leq r^{-2}h(\varphi)
+(2\pi)^{-d}\,V_F(r)\,\|\varphi\|_1^2
\label{epre1.7}
\end{equation}
for all $\varphi\in D(h)\cap L_1$ and all $r>0$.
\end{lemma}
\proof\
The proof is a slight generalization of Nash's original argument \cite{Nash}. 
We omit the details.
\hfill$\Box$

\bigskip

The  Nash inequality allows one to obtain bounds on the cross-norm $\|S_t\|_{2\to\infty}$ for many different $F$. 
In particular if $V_F$ has a polynomial behaviour one can estimate $\|S_t\|_{1\to 2}$
and  $\|S_t\|_{2\to\infty}$ by a straightforward extension of Nash's original argument.

\begin{lemma}\label{lcsg2.2}
If $h\geq f$  and $V_F(r)\leq a\,r^{(D', D)}$ then 
$\|S_t\|_{1\to\infty}\leq b\,t^{(-D/2,-D'/2)}$.
\end{lemma}
\proof\ The result can be deduced from \cite{CKS}, Theorem~2.9, or from the alternative argument given in \cite{Robm},  pages 268--269.\hfill$\Box$

\bigskip

The lemma demonstrates that  large values of $r$ give small $t$ bounds and small values of~$r$ give large $t$ bounds.
In particular  the subellipticity condition $H\geq\mu\,L^\gamma-\nu\,I$ corresponds to $ F(p)=(\mu\,|p|^{2\gamma}-\nu)\vee 0$ and this  only gives useful information on the large~$r$ behaviour of $V_F$.
It  yields bounds
$\|S_t\|_{1\to\infty}\leq a\,t^{-d/(2\gamma)}$ for $t\leq 1$.

\smallskip

Lemmas~\ref{lcsg2.1} and \ref{lcsg2.2} are  used in Section~\ref{Scsg3} to obtain uniform bounds on the semigroup kernels associated with the  Gru\v{s}in  operators  unless $n=1$ and $\delta_1$ 
or $\delta_1'$  
is in $[1/2,1\rangle$.
In the latter case one obtains subelliptic bounds of a different character, bounds in terms of the Neumann Laplacian.
But these can also be used to obtain Nash inequalities.

\smallskip

Let $F$ be a positive real function on $\Ri_+\times \Ri^m$
and let $ L_{x,N}$ denote the self-adjoint version of  the operator $-d^2/dx^2$ on $L_2(\Ri)$ with Neumann boundary conditions at the origin.
One can define the operator   $F_N=F(L_{x_1,N},i\nabla_{x_2})$ by spectral theory.
Let $f_N$ denote the corresponding quadratic form.
Now define $V_F(r)$ as the volume of the set $\{(p_1,p_2): F(p_1^2, p_2)<r^2\}$.

 \begin{lemma} \label{lcsg2.3}
If  $h\geq f_N$ in the form sense on $L_2(\Ri\times\Ri^m)$  then
\begin{equation}
\|\varphi\|_2^2\leq r^{-2}h(\varphi)
+ 4\,(2\pi)^{-d}\,V_F(r)\,\|\varphi\|_1^2
\label{epre1.8}
\end{equation}
for all $\varphi\in D(h)\cap L_1(\Ri\times\Ri^m)$ and all $r>0$.
\end{lemma}
\proof\
The proof follows the  reasoning  used to prove Lemma~\ref{lcsg2.1}
but on $L_2(\Ri_+\times \Ri^m)$.
The key point is that each $\varphi$ in the domain of the form of the Neumann Laplacian on the half-line extends by symmetry to an element in the form domain of the Laplacian on the line.
But the $L_1(\Ri)$-norm of the extension is twice the $L_1(\Ri_+)$-norm of $\varphi$.
We omit the details.
\hfill$\Box$

\section{Subelliptic estimates} \label{Scsg3}

In this section we examine the Gru\v{s}in  operators and  derive uniform estimates
on the semigroup crossnorms $\|S_t\|_{1\to\infty}$ by use of the Nash inequalities of 
Lemmas~\ref{lcsg2.1} and \ref{lcsg2.3}.
These bounds will then be improved by other techniques in Section~\ref{Scsg4}.

\begin{prop}\label{pcsg3.0}
Let $S_t$ denote the submarkovian semigroup on $L_2(\Ri^n\times\Ri^m)$ generated
by the Gru\v{s}in  operator  $H$  with coefficients $C\sim C_\delta$ where $C_\delta$  satisfies $(\ref{ecsg1.10})$.
Then
\[
\|S_t\|_{1\to\infty}\leq a\,t^{(-D/2,-D'/2)}
\]
 where
$D=(n+m(1+\delta_2-\delta_1))(1-\delta_1)^{-1}$ and 
$D'=(n+m(1+\delta_2'-\delta_1'))(1-\delta_1')^{-1}$.
\end{prop}

The local dimension $D$ depends on the parameters $\delta_1$ and $\delta_2$  governing
the local degeneracies of the coefficients of $H$ and the global dimension $D'$ depends on the parameters
$\delta_1'$ and $\delta_2'$  governing
the global  degeneracies.
Moreover, $D,D'\geq n+m$, the Euclidean dimension, with $D=n+m$ if and only if $\delta_1=0=\delta_2$
and $D'=n+m$ if and only if $\delta_1'=0=\delta_2'$.

The proof of the proposition is in two stages.
First consider the operator
\begin{equation}
H_1=-{\nabla}_x\cdot(c_{\delta_1, \delta_1'}\,\nabla_x)
\label{esub1}
\end{equation}
 on $L_2(\Ri^n)$. 
There are two possibilities.

\medskip

\noindent {\bf Either} $\delta_1\geq \delta_1'$ and then $c_{\delta_1, \delta_1'}(x)\sim 
|x|^{2\delta_1}(1+|x|^{2\delta_1})^{-1+\delta_1'/\delta_1}$,

\medskip

\noindent {\bf or} $\hspace{9mm}\delta_1\leq \delta_1'$ and then $c_{\delta_1, \delta_1'}(x)
\sim |x|^{2\delta_1}+|x|^{2\delta_1'}$.

\medskip

\noindent Since the subsequent estimates are valid up to equivalence we can effectively replace $c_{\delta_1,\delta_1'}$
by the appropriate  function on the right.
Moreover, the form
\[
h_1(\varphi)=(\nabla_x\varphi,c_{\delta_1, \delta_1'}\,\nabla_x\varphi)
=\int_{\Ri^n}dx\,c_{\delta_1, \delta_1'}(x)|(\nabla_x\varphi)(x)|^2
\]
is closed on the domain $D(h_1)=W^{1,2}(\Ri^n\,;c_{\delta_1, \delta_1'}dx)$
and $C_c^\infty(\Ri^n)$ is a core of $h$.
Therefore it suffices to establish the following form estimates on $C_c^\infty(\Ri^n)$.

\begin{prop}\label{pcsg3.1}
Let $n\geq2$, or $n=1$ and $\delta_1, \delta_1'\in[0,1/2\rangle$.
Then $h_1\geq f$ on $L_2(\Ri^n\times\Ri^m)$ where $f$ is the form of the operator $F$  with
\[
 F\sim L_x^{(1-\delta_1')}(1+ L_x)^{-(\delta_1-\delta_1')}
\]
if $\delta_1\geq\delta_1'$ and 
\[
 F\sim  L_x^{(1-\delta_1)}+  L_x^{(1-\delta_1')}
\]
if $\delta_1\leq\delta_1'$ where  $L_x=-\nabla_x^2$.

Moreover, if   $n=1$ and $\delta_1\in[1/2,1\rangle$ or $\delta_1'\in[1/2,1\rangle$ then $h_1\geq f_N$  where $f_N$ is the form of the operator $F_N$  obtained by replacing $L_x$ by $L_{x,N}$.
\end{prop}

The proof of the  subelliptic estimates of $H_1$ depends on the next two lemmas. 

\begin{lemma}\label{lsub1}
If $\gamma\in[0, 1\wedge n/2\rangle$ then
$L_x^{\,\gamma}\geq a\,|x|^{-2\gamma}$
in the form sense on $L_2(\Ri^n)$.

Moreover, if $n=1$ then $L_{x,D}\geq (4x^2)^{-1}$ on $L_2(\Ri)$ where $L_{x,D}$ is the Laplacian
with Dirichlet boundary conditions at the origin.
\end{lemma}

The statements are   versions of Hardy's inequality (see, for example, \cite{Dav13} and references
therein) and special cases of the inequalities of Caffarelli,  Kohn and Nirenberg \cite{CKN}.
The multidimensional version is often stated with $\gamma=1$ and $n\geq 3$.
In the latter case one has $a=(n-2)^2/4$ and this value is optimal. 
The fractional version follows from Strichartz' work \cite{Stri} on Fourier multipliers.

\begin{lemma}\label{lsub2}
Let $A$ and $B$ be self-adjoint operators on $L_2(\Ri^n)$ and let $\gamma\in[0, 1]$.
If  $A \geq B\geq0$ then
$A(I+A)^{-\gamma}\geq B(I+B)^{-\gamma}$
in the form sense.
\end{lemma}
\proof\
First $A\geq B$ immediately implies that 
$A(\lambda I+A)^{-1}\geq 
B(\lambda I+B)^{-1}$
for all $\lambda>0$ and this gives the result for $\gamma=1$.
But if $\gamma<1$ then
   \begin{eqnarray*}
  A(I+A)^{-\gamma}&=&c_\gamma\int^\infty_0{{d\lambda}\over{\lambda^{\gamma}}}\,A((1+\lambda)I+A)^{-1}
  \\[5pt]
  &\geq& c_\gamma\int^\infty_0{{d\lambda}\over{\lambda^{\gamma}}}\,B((1+\lambda)I+B)^{-1}
  =B(I+B)^{-\gamma}
  \end{eqnarray*}
  where we have used the standard integral representation of the fractional power.
  \hfill$\Box$
  
  \bigskip
 
 \noindent{\bf Proof of Proposition~\ref{pcsg3.1}}
Consider the case $\delta_1\geq \delta_1'$.
 Then for $\varphi\in C_c^\infty(\Ri^n)$ 
 \begin{eqnarray*}
 h_1(\varphi)&\geq&a\,(\nabla_x\varphi,|x|^{2\delta_1}(1+|x|^{2\delta_1})^{-1+\delta_1'/\delta_1}\nabla_x\varphi)\\[5pt]
 &\geq&a\,(\nabla_x\varphi,L_x^{-\delta_1}(1+L_x^{-\delta_1})^{-1+\delta_1'/\delta_1}\nabla_x\varphi)
 \geq a\,(\varphi,L_x^{1-\delta_1'}(1+L_x)^{-(\delta_1-\delta_1')}\varphi)
 \end{eqnarray*}
 by Lemmas~\ref{lsub1} and \ref{lsub2} which are applicable if $\delta_1\in[0, 1\wedge n/2\rangle$ and $\delta_1'/\delta_1\in[0,1\rangle$.
 If, however, $n=1$, $\delta_1\in[1/2,1\rangle$ and $\varphi\in C_c^\infty(\Ri^n\backslash\{0\})$  then
 \begin{eqnarray*}
 h_1(\varphi) &\geq&a\,(\varphi',L_{x,D}^{-\delta_1}(1+L_{x,D}^{-\delta_1})^{-1+\delta_1'/\delta_1}\varphi')
 \geq a\,(\varphi,L_{x,N}^{1-\delta_1'}(1+L_{x,N})^{-(\delta_1-\delta_1')}\varphi)
 \end{eqnarray*}
 where the second bound follows from the argument given in Example~5.6 of \cite{ERSZ1}.
 
 The case $\delta_1\leq \delta_1'$ is similar but simpler. 
 It also uses the basic inequality $L_x\geq L_{x,N}$ which then extends to all fractional powers.
 \hfill$\Box$

\bigskip

Next consider the   operator $H_\delta$ defined  by (\ref{ecsg1.2})
with  $\delta_1, \delta_1'\in[0,1\rangle$.
Clearly $H\geq a\, H_1$, with $H_1$ given by (\ref{esub1}) acting on $L_2(\Ri^n\times \Ri^m)$,  and so the bounds of Proposition~\ref{pcsg3.1} are applicable.
But then one has the following complementary bounds.

\begin{prop}\label{pcsg3.2}
The subelliptic estimate 
$h\geq f$ is valid on $L_2(\Ri^n\times\Ri^m)$
where $f$ is the form of the operator $F$  with
\[
 F\sim L_{x_2}^{\alpha'}(1+ L_{x_2})^{\alpha-\alpha'}
\]
where
$\alpha=(1-\delta_1)/(1+\delta_2-\delta_1)$ and
$\alpha'=(1-\delta_1')/(1+\delta_2'-\delta_1')$.
\end{prop}
\proof\
First one has $h\geq a\, h_\delta$ and after 
 a partial Fourier transformation, i.e.\  a transformation with respect to the $x_2$ variable, $H_\delta$ 
transforms to an  operator $\widetilde H_\delta=H_1+c_{\delta_2, \delta'_2}\,|p_2|^2$ on $L_2(\Ri^n)$. 
Therefore one may apply Proposition~\ref{pcsg3.1} to  $H_1$ and  bound $\widetilde H_\delta$
below by a differential operator in the $\Ri^n$-variable.
Now to proceed we again use the fractional Hardy inequality but to cover all the relevant cases
we have to pass to a fractional power of $\widetilde H_\delta$.

\begin{lemma}\label{lsub3}
Let $A$ and $B$ be positive  self-adjoint operators such that the form sum 
$A+B$ is densely defined.
Then the form sum  $A^{1/2^n}+B^{1/2^n}$ is densely defined and
\[
(A+B)^{1/2^n}\geq 2^{-1+2^{-n}}(A^{1/2^n}+B^{1/2^n})
\]
for all positive integers $n$.
\end{lemma}
\proof\
Since  $A+B\geq A$
one has  $(A+B)^{1/2}\geq A^{1/2}$ and $D(A^{1/2})\supseteq D( (A+B)^{1/2})$.
Similarly  $(A+B)^{1/2}\geq B^{1/2}$ and $D(B^{1/2})\supseteq D( (A+B)^{1/2})$.
Then
\[
A+B=(A^{1/2}+B^{1/2})^2/2+(A^{1/2}-B^{1/2})^2/2\geq (A^{1/2}+B^{1/2})^2/2
\]
and $n=1$  statement is established.
The general statement  follows by iteration.
\hfill$\Box$

\bigskip
The proof of Proposition~\ref{pcsg3.2} now continues by applying Lemma~\ref{lsub3} to deduce that
\[
\widetilde H_\delta^{1/4}
\geq  a\,(H_1^{1/4}
+c_{\delta_2, \delta'_2}^{\;\;\;\;\;1/4}\,|p_2|^{1/2})
\;.
\]
But if $\delta_1\geq\delta_1'$ and 
 $n\geq2$, or $n=1$ and $\delta_1, \delta_1'\in[0,1/2\rangle$,
then, by the first statement of Proposition~\ref{pcsg3.1},
\begin{equation}
H_1^{1/4}\geq a\,
L_{x_1}^{\,(1-\delta_1')/4}\,(I+L_{x_1}^{\,(1-\delta_1')/4})^{-1+\sigma}
\label{ecsg3.22}
\end{equation}
where $\sigma=(1-\delta_1)/(1-\delta_1')\in\langle0,1]$.
But $0\leq (1-\delta_1')/4\leq 1/4$.
Therefore 
\[
L_{x_1}^{\,(1-\delta_1')/4}\geq a\, |x_1|^{-(1-\delta_1')/2}
\]
by Lemma~\ref{lsub1}.
Hence
\begin{eqnarray*}
L_{x_1}^{\,(1-\delta_1')/4}\,(I+L_{x_1}^{\,(1-\delta_1')/4})^{-1+\sigma}
&\geq &a\,|{x_1}|^{-(1-\delta_1')/2}(1+|{x_1}|^{-(1-\delta_1')/2})^{-1+\sigma}\\[5pt]
&\geq& a\,|{x_1}|^{-(1-\delta_1)/2}
(1+|{x_1}|)^{-(\delta_1-\delta_1')/2}
\end{eqnarray*}
by Lemma~\ref{lsub2}.
Then, however,
\[
\widetilde H_\delta^{1/4}\geq a\,\Big(|{x_1}|^{-(1-\delta_1)/2}
(1+|{x_1}|)^{-(\delta_1-\delta_1')/2}
+c_{\delta_2, \delta'_2}(x_1)^{1/4}\,|p_2|^{1/2}\Big)
\;.
\]
But the right hand side is a function of $|{x_1}|$ with  a strictly positive minimum $m$  which 
is estimated by elementary arguments.
The minimum value $m$ is a positive function of $|p_2|$  which then  gives the bound $ \widetilde H_\delta\geq M  I$ with $M=m^4$.

In order to estimate the minimum $m$ we note that the first function in the last estimate, 
$x_1\mapsto |{x_1}|^{-(1-\delta_1)/2}
(1+|{x_1}|)^{-(\delta_1-\delta_1')/2}
$,  is decreasing and the second function
in the estimate, $x_1\mapsto c_{\delta_2, \delta'_2}(x_1)^{1/4}\,|p_2|^{1/2}$, is increasing.
Now if $|p_2|$ is small the graphs of the two functions intersect at a unique large value of 
$x_1\sim |p_2|^{-1/(1+\delta_2'-\delta_1')}$.
At this point the value of the sum of the functions is proportional to $|p_2|^{\alpha'/2}$.
Similarly if $|p_2|$ is large the graphs intersect at a small $x_1\sim |p_2|^{-1/(1+\delta_2-\delta_1)}$
and the minimum value is proportional to $|p_2|^{\alpha/2}$.
Therefore 
\[
m(p_2)\sim |p_2|^{\alpha'/2}(1+|p_2|^2)^{(\alpha-\alpha')/4}\;\;\;\;\;{\rm and }\;\;\;\;\;
M(p_2)\sim |p_2|^{2\alpha'}(1+|p_2|^2)^{\alpha-\alpha'}
\]
which is equivalent to the  bound stated in  the proposition.

If $n=1$ and $\delta_1\in[1/2,1\rangle$ or $\delta_1'\in[1/2,1\rangle$  the estimation procedure has to be
slightly modified.
Then Proposition~\ref{pcsg3.1}  gives the lower bound
\begin{equation}
H_1^{1/4}\geq a\,
L_{x_1,N}^{\,(1-\delta_1')/4}\,(I+L_{x_1,N}^{\,(1-\delta_1')/4})^{-1+\sigma}
\label{ecsg3.23}
\end{equation}
and $L_{x_1,N}\leq L_{x_1}$.
Nevertheless if $\alpha<1/2$ then there is an $a>0$ such that  $L_{x_1}^\alpha\geq L_{x_1,N}^\alpha\geq a\,L_{x_1}^\alpha$.
Therefore (\ref{ecsg3.22}) follows from (\ref{ecsg3.23}) by another application of 
Lemma~\ref{lsub2}.

The case $\delta_1\leq\delta_1'$ depends on the second statement of Proposition~\ref{pcsg3.1}.
The proof is similar but simpler.
\hfill$\Box$

\bigskip

Now we are prepared to estimate the crossnorm of the semigroup.

\smallskip

\noindent{\bf Proof of Proposition~\ref{pcsg3.0}}
It follows by definition that $h\sim h_\delta$ and then by  combination of Propositions~\ref{pcsg3.1} and \ref{pcsg3.2} that $h\geq f$ or $h\geq f_N$ where
$f$ is the form of a multiplier $F(L_{x_1},L_{x_2})$ and $f_N$ the form of $F(L_{x_1,N},L_{x_2})$.
Now consider the case $\delta_1\geq\delta_1'$ and assume 
 $n\geq2$, or $n=1$ and $\delta_1, \delta_1'\in[0,1/2\rangle$.
 Then one can apply Lemmas~\ref{lcsg2.1} and \ref{lcsg2.2} with 
 \[
 F(L_{x_1},L_{x_2})=a\,\Big( L_{x_1}^{(1-\delta_1')}(1+ L_{x_1})^{-(\delta_1-\delta_1')}
 +L_{x_2}^{\alpha'}(1+ L_{x_2})^{\alpha-\alpha'}\Big)
 \;.
 \]
Therefore  $V_F(r)\sim V_{F_1}(r)V_{F_2}(r)$ where 
 $V_{F_1}(r)=|\{p_1: |p_1|^{2(1-\delta_1')}(1+ |p_1|^2)^{-(\delta_1-\delta_1')}<r^2\}|$
 and
  $V_{F_2}(r)=|\{p_2: |p_2|^{2\alpha'}(1+ |p_2|^2)^{\alpha-\alpha'}<r^2\}|$.
But  $V_{F_1}(r)\sim r^{n(1/(1-\delta_1'),1/(1-\delta_1))}$.
Similarly $V_{F_2}(r)\sim r^{m(1/\alpha',1/\alpha)}$.
Therefore  $V_F(r)\sim r^{(D',D)}$.
Then the semigroup estimates of Proposition~\ref{pcsg3.0} follow from Lemma~\ref{lcsg2.2}.

The argument is similar if $\delta_1\leq\delta_1'$ but one uses the second estimate of Proposition~\ref{pcsg3.1}.
Finally if $n=1$ and $\delta_1\in[1/2,1\rangle$ or $\delta_1'\in[1/2,1\rangle$ one can make an identical argument
using the last statement of Proposition~\ref{pcsg3.1} and Lemma~\ref{lcsg2.3} which deals with the Neumann multipliers.
\hfill$\Box$


\begin{remarkn}
The arguments we have given for subellipticity estimates on $L_2(\Ri^d\,;dx)$ also extend to weighted
spaces such as  $L_2(\Ri^d\,;|x|^\beta dx)$.
In this extension the Hardy inequality is replaced by the  Caffarelli--Kohn--Nirenberg inequalities \cite{CKN}.
\end{remarkn}

\section{Comparison of kernels}\label{Scsg4}

In this section we  develop a method for transforming the uniform  bounds 
of Proposition~\ref{pcsg3.0} into bounds which better reflect the
spatial behaviour,  bounds expressed in terms of the corresponding Riemannian geometry.
In particular we establish a comparison between the Gru\v{s}in kernel and the kernel of a closely related non-degenerate operator.
Our arguments are based on wave equation techniques and are applicable to quite general
degenerate  operators.
Therefore in this section we adopt the following assumptions.

Let $C=(c_{ij})$ be   a real-valued symmetric $d\times d$-matrix with measurable, locally integrable, coefficients $c_{ij}$ which is positive-definite almost everywhere.
Again define the corresponding elliptic form $h$ by (\ref{ecsg2.1}).
The form $h$ is not necessarily closable therefore we consider its relaxation $h_0$ which is 
defined as the largest positive, closed, quadratic form $h_0$  such that $h_0\leq h$.
The relaxation occurs in the context of nonlinear phenomena and
discontinuous media  (see, for example,  \cite{Bra} \cite {ET} \cite{Jos} \cite{DalM}  \cite{Mosco}
and references therein).
It can also be defined by an approximation procedure.
 
Define $h_\varepsilon$ for each  $\varepsilon>0$ as the elliptic form corresponding to the coefficients 
$C+\varepsilon I$.
Then $h_\varepsilon$ is closed.
This follows since $h_\varepsilon$ is the  limit as $N\to\infty$ of the monotonically  increasing family of closed forms $h_{N,\varepsilon}$ 
with the bounded non-degenerate coefficients $C_{N,\varepsilon}=(C\wedge NI)+\varepsilon I$.
Moreover,   the positive selfadjoint operator $H_\varepsilon$ corresponding to the closed form $h_\varepsilon$ is the strong resolvent limit of the 
 strongly elliptic operators $H_{N, \varepsilon}$ with coefficients $C_{N,\varepsilon}$.
(For details on the limits see  \cite{BR2}, Lemma 5.2.23.)
Now the $h_\varepsilon$ are monotonically decreasing as $\varepsilon\to0$ and the 
corresponding selfadjoint operators $H_\varepsilon$ converge in the strong resolvent sense
to a positive selfadjoint operator $H_0$.
The  relaxation $h_0$ of $h$ is the Dirichlet form corresponding to $H_0$, i.e.\ 
 $D(h_0)=D(H_0^{1/2})$ and
 $h_0(\varphi)=\|H_0^{1/2}\varphi\|_2^2$.
 (For further details see the discussion in Section~2 of \cite{ERSZ1}).

 \smallskip

Next we argue that 
 the wave equation corresponding to $H_0$  has a finite speed of propagation when measured with respect to the   Riemannian (quasi-)distance  defined by
\begin{equation}
d(x\,;y)=\sup_{\psi\in D}\,(\psi(x)-\psi(y))\in[0,\infty]
\label{ecsg4.00}
\end{equation}
for all $x,y\in\Ri^d$
with
$D=\{\psi\in W^{1,\infty}(\Ri^d):\sum^d_{i,j=1}c_{ij}\,(\partial_i\psi)(\partial_j\psi)\leq1\}$.
Further, set 
\begin{equation}
d(A\,;B)=\inf_{x\in A,\,y\in B}d(x\,;y)
\label{ecsg4.100}
\end{equation}
where $A$ and $B$ are general measurable sets.

There are a variety of other methods of associating a distance with $C$ especially  if the 
coefficients are continuous. 
Then one may adopt one of several equivalent `shortest path' definitions (see \cite{JSC1} for
a survey and comparison of various possibilities for subelliptic operators).
For the current purposes the foregoing definition is most suitable.

Next we derive a Davies--Gaffney estimate \cite{Dav12} \cite{Gaf} for the semigroup generated by $H_0$  which in turn is equivalent to the  finite speed of propagation of the corresponding wave equation.
of propagation \cite{Sik} \cite{Sik3}.

\begin{prop}\label{pscsg4.0}
Let $S^{(0)}$ denote the semigroup generated by the  operator $H_0$ associated with the 
relaxation $h_0$ of the form $h$.
Then for each pair of open subsets $A,B$ of~$\Ri^d$
\begin{equation}
|(\varphi_A, S^{(0)}_t\varphi_B)|\leq e^{-d(A;B)^2(4t)^{-1}}\|\varphi_A\|_2\|\varphi_B\|_2
\label{epadty3.0}
\end{equation}
for all $\varphi_A\in L_2(A)$, $\varphi_B\in L_2(B)$ and $t>0$ with the convention
$e^{-\infty}=0$.
Moreover, the corresponding wave equation has a finite speed of propagation in the sense that
\begin{equation}
(\varphi_A,{\cos(tH_0^{1/2})}\varphi_B)=0
\label{epadty3.1}
\end{equation}
for all $\varphi_A\in L_2(A)$, $\varphi_B\in L_2(B)$ and all $t$ with  $|t|\leq d(A\,;B)$.
\end{prop}
\proof\
First let $\psi\in D$ and introduce the one-parameter family of multiplication operators
$\rho\to U_\rho=e^{\rho\psi}$ on $L_2(\Ri^d)$.
Then $\|U_\rho S^{(0)}_tU_\rho^{-1}\|_{2\to2}\leq e^{\rho^2t}$.
This is a standard estimate for strongly elliptic operators which extends to the relaxed operator
by the approximation techniques described above.
Now
\begin{eqnarray*}
|(\varphi_A,S^{(0)}_t\varphi_B)|&=&
|(U_\rho^{-1}\varphi_A,(U_\rho S^{(0)}_tU_\rho^{-1})U_\rho\varphi_B)|\\[5pt]
&\leq&e^{\rho^2t}\,\|U_\rho^{-1}\varphi_A\|_2\,\|U_\rho\varphi_B\|_2
\leq e^{-\rho d_\psi(A\,;B)}e^{\rho^2t}\,\|\varphi_A\|_2\,\|\varphi_B\|_2
\end{eqnarray*}
where
$d_\psi(A\,;B)=\inf_{x\in A,\,y\in B}(\psi(x)-\psi(y))$.
Therefore, optimizing over $\psi$  and $\rho$  one has
\begin{equation}
|(\varphi_A,S^{(0)}_t\varphi_B)|\leq  e^{-\hat d(A;B)^2(4t)^{-1}}\|\varphi_A\|_2\|\varphi_B\|_2
\label{epadty3.9}
\end{equation}
where 
\begin{equation}
\hat d(A\,;B)=\sup_{\psi\in D}d_\psi(A\,;B)
\;.
\label{epadty3.10}
\end{equation}
These estimates are valid for all measurable $A, B$ and all $\varphi_A\in L_2(A),\varphi_B\in L_2(B)$.

Since $d_\psi(A\,;B)\leq \psi(x)-\psi(y)$ for all $x\in A$ and $y\in B$ it follows that 
$\hat d(A\,;B)\leq d(A\,;B)$ again for all measurable $A$ and $B$.
But the latter inequality has a partial converse.

\begin{lemma}\label{lcsg4.0}
If  $A$ and $B$ are compact subsets then $\hat d(A\,;B)= d(A\,;B)$.
\end{lemma}
\begin{remarkns} 1. The proof is an interplay between compactness and continuity. 
It uses very little structure of the set $D$.
Indeed it suffices for the proof that $\varphi\in D$ and $c\in \Ri$ imply $-\varphi+c\in D$
and $\varphi_1,\varphi_2\in D$ imply $\varphi_1\vee \varphi_2, \varphi_1\wedge\varphi_2 \in D$.
These properties are easily verified.

\hspace{1.2cm}2. It is not necessarily the case that $\{\varphi_n\}_{n\geq1}\in D$ implies $\sup_{n\geq1} \varphi_n\in D$ or $\inf_{n\geq1}\varphi_n\in D$.
Proposition~6.5 in \cite{ERSZ1},  and its  proof, give an example of a  decreasing sequence $\chi_n$ of functions in $D$ such that $\inf_n\chi_n\not\in W^{1,\infty}$.
\end{remarkns}
\noindent{\bf Proof of Lemma~\ref{lcsg4.0}}
Since $\hat d(A\,;B)\leq d(A\,;B)$  for all measurable $A$ and $B$
it suffices to prove  $\hat d(A\,;B)\geq d(A\,;B)$ for $A$ and $B$ compact.

Fix $x\in A$ and $y\in B$.
Then for each $\varepsilon>0$ there is a $\psi_{x,y}\in D$ such that $\psi_{x,y}(x)=0$ and  
$\psi_{x,y}(y)\geq d(x\,;y)-\varepsilon/4$.
Next since $\psi_{x,y}$ is continuous there exists an open neighbourhood  $U_y$ of $y$ such that
\[
\psi_{x,y}(z)\geq d(x\,;y)-\varepsilon/2\geq d(A\,;B)-\varepsilon/2
\]
for and $z\in U_y$.
Then $B\subset \bigcup_{y\in B}U_y$ and since $B$ is compact there exist
$y_1,\ldots,y_n\in B$ such that $B\subset \bigcup_{k=1}^nU_{y_k}$.
Set $\psi_x=\sup_{1\leq k\leq n}\psi_{x,y_k}$ then $\psi_x\in D$, $\psi_x(x)=0$ and
\[
\inf_{z\in B}\psi_x(z)\geq \min_{1\leq k\leq n}d(x\,;y_k)-\varepsilon/2\geq d(A\,;B)-\varepsilon/2
\;.
\]
But since $\psi_x$ is continuous there is an open neighbourhood $U_x$ of $x$ such that $\psi_x(z)\leq \varepsilon/2$ for all $z\in U_x$.
Then by repeating the above covering  argument one can select  $x_1,\ldots,x_m\in A$ 
such that $\psi=\inf_{1\leq k\leq m}\psi_{x_k}$ satisfies $\psi(z)\leq \varepsilon/2$ for all $z\in A$.
In addition one still has $\psi(z)\geq d(A\,;B)-\varepsilon/2$ for all $z\in B$.
Therefore
\[
\inf_{x\in A,y\in B}(\psi(y)-\psi(x))\geq d(A\,;B)-\varepsilon
\;.
\]
Thus  $\hat d(A\,;B)\geq d(A\,;B)-\varepsilon$.
Since $\varepsilon>0$ is arbitrary one deduces that  $\hat d(A\,;B)\geq d(A\,;B)$ for $A$ and $B$ compact. \hfill$\Box$

\bigskip

Now we can complete the proof of Proposition~\ref{pscsg4.0}

\smallskip

\noindent{\bf End of  proof of Proposition~\ref{pscsg4.0}}
Let  $\varphi_1\in L_2(A)$ have compact support $U_1\subset A$ and $\varphi_2\in L_2(B)$
have compact support $U_2\subset B$.
Then it follows from (\ref{epadty3.9}) and Lemma~\ref{lcsg4.0} that 
\begin{eqnarray*}
|(\varphi_1,S^{(0)}_t\varphi_2)|&\leq & 
e^{-\hat d(U_1;U_2)^2(4t)^{-1}}\|\varphi_1\|_2\|\varphi_2\|_2\\[5pt]
&=&e^{- d(U_1;U_2)^2(4t)^{-1}}\|\varphi_1\|_2\|\varphi_2\|_2
\leq e^{- d(A;B)^2(4t)^{-1}}\|\varphi_1\|_2\|\varphi_2\|_2
\;.
\end{eqnarray*}
But if $A$ is  open  then the functions of compact support in $L_2(A)$ are dense
and similarly for $B$.
Therefore the first statement  (\ref{epadty3.0}) of the proposition follows by continuity.
The second statement (\ref{epadty3.1})  is a direct consequence (see, for example, \cite{ERSZ1} Lemma~3.3).
\hfill$\Box$

\bigskip

An equivalent way of expressing the finite speed of  propagation  is the following.
\begin{lemma} \label{lpre2.11}
Let $A$  be an  open subset and $F$ a closed subset   with $A\subset F$.
Then 
\[
{\cos(tH_0^{1/2})}L_2(A)\subseteq L_2(F)
\]
for all $t\in\Ri$ with $|t|\leq d(A\,;F^{\rm c})$.
\end{lemma}
\proof\
Let $B$ be an  open subset of $ F^{\rm c}$. Then $d(A\,;F^{\rm c})\leq d(A\,;B)$.
Therefore $L_2(B)\perp {\cos(tH_0^{1/2})}L_2(A)$ for all $t\in\Ri$ with $|t|\leq d(A\,;F^{\rm c})$
by (\ref{epadty3.1}).
Hence $L_2(F^{\rm c})\perp {\cos(tH_0^{1/2})}L_2(A)$ and one must have 
${\cos(tH_0^{1/2})}L_2(A)\subseteq L_2(F)$.\hfill$\Box$

\bigskip

  The subsequent comparison theorem depends on a generalization of the propagation property which
  emphasizes the local nature.
  As a preliminary let $h_1$ and $h_2$ be two  elliptic forms with coefficients $C_1=(c^{(1)}_{ij})$ and $C_2=(c^{(2)}_{ij})$ and let $H_{1,0}$ and $H_{2,0}$ denote the corresponding relaxations.
  Moreover, assume that $C_1\geq C_2$ and $C_1\geq\mu I>0$.
  In particular $D(h_1)\subseteq  D(h_2)\subseteq W^{1,2}(\Ri^d)$ and $d_1(x\,;y)\leq d_2(x\,;y)$ for all $x,y$ where $d_1$ and $d_2$ denote the Riemannian distances associated with $C_1$ and $C_2$.
  Let $U=\supp (C_1-C_2)$ and set $F=\overline{U^{\rm c}}$.
  
  \begin{lemma}\label{lcsg4.2} If $A$ is an open subset of the closed subset $F$ then
  \begin{equation}
\cos(tH_{1,0}^{1/2})\varphi_A=\cos(tH_{2,0}^{1/2})\varphi_A
\label{ecsg4.5}
\end{equation}
for all $\varphi_A\in L_2(A)$ and all $t\in \Ri$ with $|t|\leq d_1(A\,;U)$.
\end{lemma}
\proof\
If $C_1$ and $C_2$ are strongly elliptic this result follows from the proof  of Proposition~3.15 in 
\cite{ERS4}.
The extension to the more general situation can then be made by approximation as follows.

Let $N>\varepsilon>0$.
Set $C_{1,N,\varepsilon}=(C_1\wedge N I)+\varepsilon I$ and $C_{2,N,\varepsilon}=(C_2\wedge N I)+\varepsilon I$.
Then the lemma is valid for the corresponding strongly elliptic operators $H_{1,N,\varepsilon}$,
$H_{2,N,\varepsilon}$.
But as $N\to\infty$ these operators converge in the strong resolvent sense to 
$H_{1,\varepsilon}=H_1+\varepsilon I$ and $H_{2,\varepsilon}=H_2+\varepsilon I$, respectively.
Since $H_{1,\varepsilon}\geq H_{1,N,\varepsilon}\geq H_{2,N,\varepsilon}$  it follows that 
\[
\cos(tH_{1,\varepsilon}^{1/2})\varphi_A=\cos(tH_{2,\varepsilon}^{1/2})\varphi_A
\]
for all $t\in\Ri$ with $|t|\leq d_{1,\varepsilon}(A\,;U)$ where $d_{1,\varepsilon}$
denotes the Riemannian distance associated with $C_{1,\varepsilon}$.

Finally it follows that $H_{1,\varepsilon}$ and $H_{2,\varepsilon}$ converge in the strong resolvent sense to $H_{1,0}$ and $H_{2,0}$, respectively, as $\varepsilon\to0$.
Moreover, since $C_1\leq C_{1,\varepsilon}\leq (1+\varepsilon\mu^{-1}) \,C_1$ it follows that
$(1+\varepsilon\mu^{-1})^{-1/2}d_1(A\,;B)\leq d_{1,\varepsilon}(A\,;B)\leq d_1(A\,;B)$
for all measurable $A$, $B$.
Hence $d_{1,\varepsilon}(A\,;B)\to d_1(A\,;B)$ as $\varepsilon\to0$.
The statement of the lemma follows in the limit.
\hfill$\Box$

  \bigskip

Now we are prepared to establish the principal comparison result.
In the sequel $S^{(1,0)}$, $S^{(2,0)}$ denote the semigroups generated by 
$H_{1,0}$, $H_{2,0}$  and  $K^{(1,0)}$, $K^{(2,0)}$ denote the corresponding
kernels.

\begin{thm}\label{tcsg3.1}
Adopt the foregoing notation and assumptions.
Let 
$\one_A$  denote the indicator function of the open subset 
$A$ of $F=\overline{U^{\rm c}}$ where $U=\supp (C_1-C_2)$.
  Set 
\[
M_N(t)=\|\one_A(I+t^2\,H_{1,0})^{-N}\one_A\|_{1\to\infty}+\|\one_A(I+t^2\,H_{2,0})^{-N}\one_A\|_{1\to\infty}
\]
for $N\in\Ni$.

Then there is an $a_N>0$  such that 
\begin{equation}
\sup_{x,y \in A} |K^{(1,0)}_t(x\,;y)- K^{(2,0)}_t(x\,;y)| \leq a_N\,
M_N(t/\rho )\,(\rho^2/t)^{-1/2}
 \,e^{-\rho^2/(4t)}
\label{ecsg3.2}
\end{equation}
for all $t>0$ where $\rho =d_1(A\,;U)$.
\end{thm}

Note that there is no reason that $M_N$ is finite.
Subsequently we give conditions which
ensure that $M_N$ is indeed finite for large $N$.

The proof of the theorem is based on the estimates developed in \cite{Sik} and \cite{Sik3}.
As a preliminary we need some properties of functions of the operators $H_{0, 1}$ and 
$H_{2,0}$.

First, let $\Psi$ be an  even bounded Borel function with
Fourier transform 
$\widetilde \Psi$ satisfying $\supp \widetilde\Psi\subseteq [-\rho, \rho]$ where $\rho>0$.
Then for each pair of  open subsets $B_1$ and $B_2$ 
\begin{equation}
(\varphi_{B_1},\Psi(H_{1,0}^{1/2})\varphi_{B_2})=0=(\varphi_{B_1},\Psi(H_{2,0}^{1/2})\varphi_{B_2})
\label{ecsg3.4}
\end{equation}
for all $\varphi_{B_1}\in L_2(B_1)$, $\varphi_{B_2}\in L_2(B_2)$ where  
$\rho\leq d_1(B_1\,;B_2)$.
This follows for $H_{1,0}$ from the representation
\[
 \Psi(H_{1,0}^{1/2})
=(2\pi)^{-1/2}\int_{\Ri}dt\, \widetilde \Psi(t) \exp(itH_{1,0}^{1/2})
=(2\pi)^{-1/2}\int_{-\rho}^\rho dt\, \widetilde \Psi(t) \cos(tH_{1,0}^{1/2})\;.
\]
and   condition (\ref{epadty3.1}).
The argument for $H_{2,0}$ is similar.

\begin{lemma}\label{lcsg3.1}
Let $\Psi$ be an  even bounded Borel function with
Fourier transform 
$\widetilde \Psi$ satisfying $\supp \widetilde\Psi\subseteq [-2\rho ,2\rho ]$.
Then
\[
\one_A \Psi(H_{1,0}^{1/2})\one_A=\one_A \Psi(H_{2,0}^{1/2})\one_A
\;.
\]
\end{lemma}
\proof\
Since  $A\subset F$
and  $\rho =d_1(A\,;U)$ one has 
\begin{equation}
\cos(tH_{1,0}^{1/2})\varphi_A=\cos(tH_{2,0}^{1/2})\varphi_A
\label{ecsg3.5}
\end{equation}
for all  $\varphi_A\in L_2(A)$ and $t\in \Ri$ with $|t|\leq \rho$ by Lemma~\ref{lcsg4.2}.

Next remark that 
\begin{eqnarray*}
(\psi_A,\cos(2tH_{1,0}^{1/2})\varphi_A)-(\psi_A,\cos(2tH_{2,0}^{1/2})\varphi_A)\\[5pt]
&&\hspace{-5cm}=
2\,\Big((\cos(tH_{1,0}^{1/2})\psi_A,\cos(tH_{1,0}^{1/2})\varphi_A)-
(\cos(tH_{2,0}^{1/2})\psi_A,\cos(tH_{2,0}^{1/2})\varphi_A)\Big)=0
\end{eqnarray*}
for $|t|\leq \rho$ and $\psi_A,\varphi_A\in L_2(A)$.
Then, however, one has
\begin{eqnarray*}
(\psi, \one_A \Psi(H_{1,0}^{1/2})\one_A\varphi)
&=&(2\pi)^{-1/2}\int_{-2\rho}^{2\rho}dt\, \widetilde \Psi(t) (\one_A\psi, \cos(tH_{1,0}^{1/2})\one_A\varphi)\\[5pt]
&=&(2\pi)^{-1/2}\int_{-2\rho}^{2\rho}dt\, \widetilde \Psi(t) (\one_A\psi,\cos(tH_{2,0}^{1/2})\one_A\varphi)=
(\psi, \one_A \Psi(H_{2,0}^{1/2})\one_A\varphi)
\end{eqnarray*}
for all $\varphi\in L_2(\Ri^d)$.
\hfill$\Box$

\bigskip

Now we are prepared to prove the theorem.
We use the notation $K_T$ for the kernel of an operator $T$.

\smallskip

 \noindent{\bf Proof of Theorem~\ref{tcsg3.1}}
 Let  $\psi \in C^{\infty}(\Ri)$ be an increasing function with
\[
\psi(u) = \left\{ \begin{array}{ll}
    0 & \mbox{ if $u \le -1$}\\
    1  & \mbox{ if $u \ge -1/2$}\; .
           \end{array}
    \right.
\]
Then for $s>1$ define the family of functions $\varphi_s$ such that 
$\varphi_s(u)= \psi(s(|u|-s))$.
Next  define functions $\widetilde\Phi_s$ and $\widetilde\Psi_s$ by 
\begin{equation}
\widetilde\Phi_s(u)
= (4\pi)^{-1/2}\exp {({-u^2/4})}-\widetilde\Psi_s(u)
= \varphi_s(u) \,(4\pi)^{-1/2}\exp {({-u^2/4})}
\;.
\end{equation}
Then  the inverse  Fourier transforms satisfy 
${\Phi_s}(\lambda)+{\Psi_s}(\lambda)=(2\pi)^{-1/2}\,\exp(-\lambda^2)$
and
\begin{equation}
S^{(i,0)}_t=\exp(-tH_{i,0})={\Phi_s}((tH_{i,0})^{1/2})+{\Psi_s}((tH_{i,0})^{1/2})
\label{ggg}
\end{equation}
for $i=1,2$.
Integration by parts $2N$ times yields
\begin{eqnarray*}
\int du\,e^{-u^2/4}e^{-iu\lambda}\,\varphi_s(u)=
\int du\,e^{-u^2/4-i\lambda u}\,\underbrace{\Big(\frac{1}{u/2+i\lambda}%
\Big(\ldots\Big(\frac{1}{u/2+i\lambda}\varphi_s(u)\Big)^{\prime}\ldots\Big)'\Big)%
^{\prime}}_{2N}\;.
\end{eqnarray*}
Hence for any $N\in\Ni$ and $s>1$ there is an $a_N>0$ such that
\begin{equation}\label{osz}
|{\Phi_s}(\lambda)| \le
a_N\,\frac{1}{s\,(1+\lambda^2/s^2)^{N}}\,e^{-s^2/4} \;\;\;,\label{osz1}
\end{equation}
with the value of $a_N$ depending only on $N$.

Next note that
supp~$\widetilde\Psi_s \subseteq [-s+(2s)^{-1},s-(2s)^{-1}]\subset [-s,s]$.
So setting $s_{\rho}= 2\,\rho\,t^{-1/2} $ one has 
\[
\one_A{{\Psi_{s_{\rho}}}((tH_{1,0})^{1/2})}\one_A=
\one_A{{\Psi_{s_{\rho}}}((tH_{2,0})^{1/2})}\one_A
\]
by Lemma~\ref{lcsg3.1} and rescaling with $t^{1/2}$.
Hence  one deduces from (\ref{ggg}) that 
\begin{eqnarray}
\one_A(S^{(1,0)}_t-S^{(2,0)}_t)\one_A&=&
\one_A{{\Phi_{s_{\rho}}}((tH_{1,0})^{1/2}))}\one_A-\one_A{{\Phi_{s_{\rho}}}((tH_{2,0})^{1/2}))}\one_A
\;.
\label{incl}
\end{eqnarray}
Therefore
\begin{eqnarray*}
\sup_{x,y \in A} |K^{(1,0)}_t(x\,;y)- K^{(2,0)}_t(x\,;y)|&&\\[5pt]
& &\hspace{-2cm}{}\leq
\|\one_A{\Phi_{s_{\rho}}}((tH_{1,0})^{1/2})\one_A\|_{1\to\infty}+
\|\one_A{\Phi_{s_{\rho}}}((tH_{2,0})^{1/2})\one_A\|_{1\to\infty}
\;.
\end{eqnarray*}

Now let $\Theta_{s_{\rho}}$ be a possibly complex function such that
$\Theta_{s_{\rho}}(\lambda)^2={\Phi_{s_{\rho}}}(t^{1/2}\lambda)$.
Then
\begin{eqnarray*}
\|\one_A{\Phi_{s_{\rho}}}((tH_{i,0})^{1/2}))\one_A\|_{1\to\infty}&=&
(\|\one_A\,\Theta_{s_{\rho}}((H_{i,0})^{1/2})\|_{2\to\infty})^2\\[15pt]
&&\hspace{-3cm}\leq (\|(I+t^2H_{i,0}/\rho^2)^{N/2}\,\Theta_{s_{\rho}}((H_{i,0})^{1/2})\|_{2\to2})^2\,(\|\one_A(I+t^2H_{i,0}/\rho^2)^{-N/2}\|_{2\to\infty})^2
\;.
\end{eqnarray*}
But
\begin{eqnarray*}
\left|\Theta_{s_{\rho}}(\lambda)
(1+t^2\lambda^2/\rho^2)^{N/2}\right|^2
&=&
\left|{\Phi_{s_{\rho}}}(t^{1/2}\lambda)
(1+t(t^{1/2}\lambda)^2/\rho^2)^{N}\right|\\[5pt]
&=&
\left|{\Phi_{s_{\rho}}}(t^{1/2}\lambda)
(1+4(t^{1/2}\lambda)^2/{s_\rho}^2)^{N}\right|\;.
\end{eqnarray*}
Therefore
\begin{eqnarray*}
(\|(I+t^2H_{i,0}/\rho^2)^{N/2}\,\Theta_{s_{\rho}}((H_{i,0})^{1/2})\|_{2\to2})^2&=&
\sup_{\lambda \ge 0}\left|\Theta_{s_{\rho}}(\lambda)
(1+t^2\lambda^2/\rho^2)^{N/2}\right|^2\\[5pt]
&=&\sup_{\lambda \ge 0}
\left|{\Phi_{s_{\rho}}}(t^{1/2}\lambda)
(1+4(t^{1/2}\lambda)^2/{s_\rho}^2)^{N}\right|\\[5pt]
&\leq& a_N\,(s_\rho/2)^{-1}\,e^{-s_\rho^2/16}=a_N\,(\rho^2/t)^{-1/2}\,e^{-\rho^2/(4t)}
\;.
\end{eqnarray*}
Combining these estimates gives
\begin{eqnarray*}
\sup_{x,y \in A} |K^{(1,0)}_t(x\,;y)- K^{(2,0)}_t(x\,;y)|
&=&a_N\,M_N(t/\rho)\,(\rho^2/t)^{-1/2}\,e^{-\rho^2/(4t)}
\end{eqnarray*}
which establishes the statement of the theorem.
\hfill$\Box$

\bigskip

The statement of the theorem can be reformulated in terms of {\it a priori}
bounds on the semigroups $S^{(1,0)}$ and $S^{(2,0)}$ if one has 
 suitable uniform estimates on the crossnorms $\|S^{(i,0)}\|_{1\to\infty}$.

 If $h$ is strongly elliptic then  $\|S^{(0)}_t\|_{1\to\infty}\leq a\,t^{-d/2}$ for all $t>0$ and we will
 assume  analogous bounds 
$ \|S^{(0)}_t\|_{1\to\infty}\leq a\,V(t)^{-1}$
  where $V$ is a positive increasing function which satisfies the usual doubling property
 \begin{equation}
V(2t)\leq a\,V(t)
\label{epre1.1}
\end{equation}
for some $a>0$ and all $t>0$.
It  follows from (\ref{epre1.1})  that there are $a, \widetilde{D}>0$  such that 
\begin{equation}
V(s)\leq a \,(s/t)^{\tilde{D}}\,V(t)
\label{epre1.17}
\end{equation}
for all $s\geq t>0$.

The parameter $\widetilde{D}$ is referred to as the doubling dimension, although it need not be an integer.

\begin{cor}\label{ccsg3.1}
Adopt the hypotheses and notation of Theorem~$\ref{tcsg3.1}$.
Let $V$ be a positive increasing function which satisfies the doubling property~$(\ref{epre1.1})$.
Assume that 
\[
\|\one_A\,S^{(1,0)}_t\,\one_A\|_{1\to\infty}\vee \|\one_A\,S^{(2,0)}_t\,\one_A\|_{1\to\infty}\leq V(t)^{-1}
\]
for all $t>0$.

Then there is an $a>0$ such that 
\begin{equation}
\sup_{x,y \in A} |K^{(1,0)}_t(x\,;y)- K^{(2,0)}_t(x\,;y)| \leq a\,
V(t^2/\rho^2)^{-1}\,(\rho^2/t)^{-1/2}
 \,e^{-\rho^2/(4t)}
\label{ecsg3.20}
\end{equation}
for all $t>0$, where $\rho=d_1(A\,;U)$.
\end{cor}
\proof\
It follows by  Laplace transformation that
\begin{eqnarray*}
\|\one_A (I+t^2\,H_0)^{-N}\one_A\|_{1\to\infty}&\leq& 
{{1}\over{(N-1)!}}\int^\infty_0ds\,s^{N-1}e^{-s}\|\one_A\,
 S^{(0)}_{st^2}\,\one_A\|_{1\to\infty}\\[5pt]
&\leq &{{1}\over{(N-1)!}}\int^\infty_0ds\,s^{N-1}e^{-s}\,V(st^2)^{-1}\leq a_N\,V(t^2)^{-1}
\end{eqnarray*}
for all $t>0$ with $a_N$ finite if $N>\tilde  D$.
The last step uses the doubling property (\ref{epre1.17}) in the form $V(t^2)\leq \rho\,s^{-\tilde{D}}\,V(st^2)$ for $s\leq 1$.
Therefore $M_N(t)\leq 2\,a_N\,V(t^2)^{-1}$ and the statement of the corollary is an immediate
consequence of Theorem~\ref{tcsg3.1}.
\hfill$\Box$

\bigskip

In Section~\ref{Scsg6} we apply Corollary~\ref{ccsg3.1} to the Gru\v{s}in operator $H$.
Then we set $H_2=H$ and $H_1$ is a Gru\v{s}in operator  with no local degeneracies
but with the same growth properties for $|x_1|\geq 1$.

\section{Volume estimates}\label{Scsg5}

In this section we  estimate  the Riemannian distance  $d(\cdot\,;\cdot)$
associated with a general Gru\v{s}in
operator $H$  and   the  volume of the   balls
$B(x\,;r)=\{y\in\Ri^n\times\Ri^m: d(x\,;y)<r\}$.
In particular we prove that the balls have the volume doubling property.

First, remark that if $C_1$, $C_2$ are two positive symmetric matrices whose entries
are measurable functions and $d_1(\cdot\,;\cdot)$, $d_2(\cdot\,;\cdot)$ the corresponding distances
then $C_1\sim C_2$ implies  $d_1(\cdot\,;\cdot)\sim d_2(\cdot\,;\cdot)$.
In particular if $a\,C_1\leq C_2\leq b\,C_1$ then $b^{-1/2}d_1(x\,;y)\leq d_2(x\,;y)\leq a^{-1/2}d_1(x\,;y)$
for all $x,y\in \Ri^d$.
Moreover, the corresponding balls $B_1$, $B_2$ satisfy
\[
B_2(x\,;b^{-1/2}r)\subseteq B_1(x\,;r)\subseteq B_2(x\,;a^{-1/2}r)
\]
for all $x\in \Ri^d$ and all $r>0$.
Then $|B_1|\sim |B_2|$ and if $|B_1|$ satisfies the doubling property with doubling  dimension $\tilde D$ then $|B_2|$
also satisfies the property with the same dimension.

Secondly, the coefficient matrix $C$ of the Gru\v{s}in operator  satisfies $C\sim C_\delta$
and to calculate the Riemannian distance, up to equivalence, we may make a convenient choice of the
$C_\delta$.
In particular we may choose $C_\delta$ such that its entries $c_{\delta_1,\delta'_1}$,  $c_{\delta_2,\delta'_2}$ are continuous functions over $\Ri^n$.
In fact we may assume $c_{\delta_1,\delta'_1}\in C^{\delta_1}(\Ri^n)$ and 
 $c_{\delta_2,\delta'_2}\in C^{\delta_2}(\Ri^n)$.
 The continuity of the coefficients then allow us to appeal to path arguments
 in the computation of the Riemannian distance.
 Specifically the Riemannian distance  (\ref{ecsg4.00}) is equivalent to the shortest distance
 of paths measured with respect to a continuous choice of $C_\delta$
 (see, for example, \cite{JSC1}).
 The restriction $\delta_1\in [0,1\rangle$ is essential to ensure that there is a continuous path between
 each pair of points.

Thirdly, let  $B(x_1,x_2\,;r)$  denote the ball with centre $x=(x_1,x_2)$
and    $|B(x_1,x_2\,;r)|$  its volume.
Further let  $D$
and  $D'$
denote the parameters occurring in the uniform bounds of
Proposition~\ref{pcsg3.0}.
Next define the function   $\Delta_\delta$ by the formula
\begin{equation}\label{csgdelta}
  \Delta_\delta (x_1,x_2\,;y_1,y_2)=   \left\{ \begin{array}{llll}
{|x_2-y_2|}/{(|x_1|+|y_1|)^{(\delta_2,\delta_2')}}
 & \mbox{ if $|x_2-y_2| \leq  (|x_1|+|y_1|)^{(\rho,\rho')}  $}\\[8pt]
  {|x_2-y_2|^{(1-\gamma,1-\gamma')}  }
&     \mbox{ if
$ |x_2-y_2|\geq (|x_1|+|y_1|) ^{(\rho,\rho')}   $}  
           \end{array}
    \right.
\end{equation}
where $\rho=1+\delta_2-\delta_1$,  $\rho'=1+\delta_2'-\delta_1'$,
$\gamma =\delta_2/\rho$ and 
$\gamma' =\delta'_2/\rho'$.
Now  set
$$
  D_\delta (x_1,x_2\,;y_1,y_2)=
{|x_1-y_1|}/{(|x_1|+|y_1|)^{(\delta_1,\delta_1')}}
+\Delta_\delta (x_1,x_2\,;y_1,y_2).
$$
Note that
$$
 \Delta_\delta (x_1,x_2\,;y_1,y_2)\sim
\frac{|x_2-y_2|}{ (|x_1|+|y_1|)^{(\delta_2,\delta_2')}+
|x_2-y_2|^{(\gamma,\gamma')} }.
$$
Note also that if  $ (|x_1|+|y_1|) ^{(\rho,\rho')}= |x_2-y_2|   $
then
$$
\frac{|x_2-y_2|^{\phantom{(\delta_2,\delta_2')}}}{(|x_1|+|y_1|)^{(\delta_2,\delta_2')}}
=\frac{|x_2-y_2|^{\phantom{(\gamma,\gamma')}}}{|x_2-y_2|^{(\gamma,\gamma')}
}=|x_2-y_2|^{(1-\gamma,1-\gamma')}
$$
so $ \Delta_\delta$ is a continuous function of  the variables
$x_1,x_2,y_1,y_2$.
Now we can make a very explicit estimate of the Riemannian distance associated with the coefficients $C_\delta$.
Estimates of the same general nature have been given for a different class of Gru\v{s}in operators
in \cite{FGW2}, Proposition 2.2 (see also \cite{Fra}, Theorem~2.3).

\begin{prop}\label{pcsg5.1}
Consider the general Gru\v{s}in operator with coefficients $C\sim C_\delta$.
If $d_\delta$ is the Riemannian distance 
$ (\ref{ecsg4.00})$ but with coefficients $C_\delta$  then
\[
d_\delta(x_1,x_2\,;y_1,y_2)\sim D_\delta (x_1,x_2\,;y_1,y_2)
\;.
\]
Moreover, the volume of the corresponding  balls satisfy
\begin{equation}\label{evol}
|B(x_1,x_2;r)|\sim \left\{ \begin{array}{llll}
 r^{(D,D')} & \mbox{ if $r\geq |x_1|^{(1-\delta_1,1-\delta_1')} $}\\[5pt]
   r^{n+m}|x_1|^{(\beta,\beta')}
& \mbox{ if $ r \leq |x_1|^{(1-\delta_1,1-\delta_1')} $}
           \end{array}
    \right.
\end{equation}
where
$\beta=n\delta_1+m\delta_2$ and $\beta'=n\delta'_1+m\delta_2'$.
\end{prop}
\proof\
First note that the coefficients of  $C_\delta$ do not depend on $x_2$. 
Hence
\[
d_\delta(x_1,x_2\,;y_1,y_2)=d_\delta(x_1,0\,;y_1,y_2-x_2)
\;.
\]
Without loss of generality one may  assume  $|x_1| \leq |y_1|$
and so $|y_1|\sim |x_1|+|y_1|$. 
We adopt  this convention throughout the remainder of the  proof.
Next by the  triangle inequality
\begin{equation}
d_\delta(x_1,0\,;y_1,y_2-x_2)\leq
d_\delta(x_1,0\,;y_1,0)+ d_\delta(y_1,0\,;y_1,y_2-x_2)
\;.\label{egr5.1}
\end{equation}
Now we  argue that the first term on the right hand side of (\ref{egr5.1})  satisfies the estimate
\[
d_\delta(x_1,0\,;y_1,0)\leq a\,
\frac{|x_1-y_1|^{\phantom{(\delta_1,\delta_1')}}}
{(|x_1|+|y_1|)^{(\delta_1,\delta_1')}}
\;.
\]
In order to establish this inequality we distinguish between two cases:
 $|x_1-y_1| \geq |y_1|/2$ and  $|x_1-y_1| \leq |y_1|/2$.

If $|x_1-y_1| \ge |y_1|/2$    then
$d_\delta(x_1,0\,;y_1,0)\leq 2\, d_\delta(0,0\,;y_1,0)$.
But $d_\delta(0,0\,;y_1,0)$ is less than the length of a straight line path
from $0$ to $y_1$.
Thus 
\[
d_\delta(x_1,0\,;y_1,0)\leq a\,
\frac{|y_1|^{\phantom{(\delta_1,\delta_1')}}}{|y_1|^{(\delta_1,\delta_1')}  }\leq
a\,\frac{|x_1-y_1|^{\phantom{(\delta_1,\delta_1')} }}{(|x_1|+|y_1|)^{(\delta_1,\delta_1')}  }
\;.
\]
If, however,   $|x_1-y_1| \leq |y_1|/2$ then    we consider the path
$(x_1(t),x_2(t))=(tx_1+(1-t)y_1,0)$. 
Note that
 $|tx_1+(1-t)y_1| \le |y_1|-t |x_1-y_1|  \leq  |y_2|/2$
so
\[
d_\delta(x_1,0\,;y_1,0)\leq a \int_0^1dt\, |x_1-y_1| \,c_{\delta_1,\delta_1'}(|y_1|/2)^{-1/2}
\sim \frac{|x_1-y_1|^{\phantom{(\delta_1,\delta_1')} }}{(|x_1|+|y_1|)^{(\delta_1,\delta_1')}  }
\;.
\]
This completes the bound of the first term on the right hand side of (\ref{egr5.1}).

Next we bound the second term on the right  of (\ref{egr5.1}).
Specifically we will establish that
\begin{equation}
 d_\delta(y_1,0\,;y_1,y_2-x_2) \le a\, \Delta_\delta (y_1,0\,;y_1,x_2-y_2)\sim
\Delta_\delta (x_1,x_2\,;y_1,y_2)
\;.\label{csg.one}
\end{equation}
If $|x_2-y_2|\le |y_1|^{(\rho,\rho')} $
then considering the path $y(t)=(y_1,t(x_2-y_2))$ we find
\begin{eqnarray}
 d_\delta(y_1,0\,;y_1,y_2-x_2)&\leq&
a\int_0^1dt\,  |x_2-y_2|\,c_{\delta_2,\delta_2'}(|y_1|)^{-1/2}\nonumber\\[5pt]
&\sim &{|x_2-y_2|}{|y_1|^{(-\delta_2,-\delta_2')} }
\sim \Delta_\delta (x_1,x_2\,;y_1,y_2)
\;.
\label{ecsg.a1}
\end{eqnarray}
If, however,
$|x_2-y_2|\geq |y_1|^{(\rho,\rho')} $ we set $\tilde{y}_1=
(y_1/{|y_1|})|x_2-y_2|^{(1/\rho,1/\rho')}$.
Then by (\ref{ecsg.a1})
\[
d_\delta(\tilde{y}_1,0\,;\tilde{y}_1, y_2-x_2)\leq
a\,{|x_2-y_2|}|\tilde{y}_1|^{(-\delta_2,-\delta_2')}  \sim
{|x_2-y_2|^{(1-\gamma,1-\gamma')}}
\;.
\]
Therefore
\begin{eqnarray*}
 d_\delta(y_1,0\,;y_1,y_2-x_2)&\leq&
2 \,d_\delta(\tilde{y}_1,0;\tilde{y}_1,0)+
d_\delta(\tilde{y}_1,0;\tilde{y}_1, y_2-x_2)\\[5pt]
&\leq& a\, ( |\tilde{y}_1|^{(1-\delta_1,1-\delta_1')}+
{|x_2-y_2|^{(1-\gamma,1-\gamma')}})\\[5pt]
&\sim& {|x_2-y_2|^{(1-\gamma,1-\gamma')}}
\sim  \Delta_\delta (x_1,x_2\,;y_1,y_2)
\;.
\end{eqnarray*}
Combination of these estimates then gives an upper bound
\[
d_\delta(x_1,x_2\,;y_1,y_2)\leq a\, D_\delta (x_1,x_2\,;y_1,y_2)
\]
for all $x,y$.
Therefore to complete the proof of equivalence of the distances
we have to establish a similar lower bound.

\smallskip

First we argue that
\begin{equation}
  c_{\delta_1,\delta_1'}(x_1)(\nabla_{x_1}D_\delta (x_1,x_2;y_1,y_2))^2
+c_{\delta_2,\delta_2'}(x_1)(\nabla_{x_2}D_\delta (x_1,x_2;y_1,y_2))^2
\leq a
\;.\label{csg.nabla}
\end{equation}
To establish (\ref{csg.nabla}) we first note that
$\nabla_{x_2}({|x_1-y_1|}/{(|x_1|+|y_1|)^{(\delta_1,\delta_1')}})=0$.
Moreover,
\begin{eqnarray*}
\nabla_{x_1}\frac{|x_1-y_1|^{\phantom{(\delta_1,\delta_1')}}}
{(|x_1|+|y_1|)^{(\delta_1,\delta_1')}}=
({\nabla_{x_1}|x_1-y_1|})\,{(|x_1|+|y_1|)^{ (-\delta_1,-\delta_1') }}+
|x_1-y_1|\,\nabla_{x_1}(|x_1|+|y_1|)^{(-\delta_1,-\delta_1')}\\[5pt]
 \le  (|x_1|+|y_1|)^{(-\delta_1,-\delta_1') }  +a\,
\frac{|x_1-y_1|^{\phantom{(1+\delta_1,1+\delta_1')}}}
{(|x_1|+|y_1|)^{(1+\delta_1,1+\delta_1')}}
\le a\,(|x_1|+|y_1|)^{(-\delta_1,-\delta_1')}\;.
\end{eqnarray*}
Therefore
\[
 c_{\delta_1,\delta_1'}(x_1)\left(\nabla_{x_1}
\frac{|x_1-y_1|^{\phantom{(\delta_1,\delta_1')}}}{(|x_1|+|y_1|)^{(\delta_1,\delta_1')}}\right)^2
+c_{\delta_2,\delta_2'}(x_1)\left(\nabla_{x_2}
\frac{|x_1-y_1|^{\phantom{(\delta_1,\delta_1')}}}{(|x_1|+|y_1|)^{(\delta_1,\delta_1')}}\right)^2 \le a
\;.
\]
Hence to establish  (\ref{csg.nabla}) it is enough to show that
\begin{equation}\label{csg.nabla1}
  c_{\delta_1,\delta_1'}(x_1)(\nabla_{x_1}\Delta_\delta
(x_1,x_2;y_1,y_2))^2
+c_{\delta_2,\delta_2'}(x_1)(\nabla_{x_2}\Delta_\delta
(x_1,x_2;y_1,y_2))^2
\le a
\;.
\end{equation}
Now
$c_{\delta_1,\delta_1'}(x_1)(\nabla_{x_1}|x_2-y_2|^{(1-\gamma,1-\gamma')})^2=0$
and
\[
c_{\delta_2,\delta_2'}(x_1)(\nabla_{x_2}|x_2-y_2|^{(1-\gamma,1-\gamma')})^2
\leq \left(\frac{|x_1|^{(\delta_2,\delta_2')}}{
|x_2-y_2|^{(\gamma,\gamma')}}\right)^2
\le a
\]
for all  $|x_2-y_2|\ge |x_1|+|y_1|^{(\rho,\rho')} $.
Next
\begin{eqnarray*}
   c_{\delta_1,\delta_1'}(x_1) \left(\nabla_{x_1}
\frac{|x_2-y_2|^{\phantom{(\delta_2,\delta_2')}} }{  (|x_1|+|y_1|)^{(\delta_2,\delta_2')}   }\right)^2
&\leq& a\,  c_{\delta_1,\delta_1'}(x_1)
\left(\frac{|x_2-y_2|^{ \phantom{(1+\delta_2,1+\delta_2')} }}{(|x_1|+|y_1|)^{ (1+\delta_2,1+\delta_2') }
}\right)^2\\[5pt]
&&\hspace{-2cm}{}\leq a\,
\left(\frac{|x_2-y_2|\,|x_1|^{ (\delta_1,\delta_1')\phantom{1} } }{ (|x_1|+|y_1|)^{
(1+\delta_2,1+\delta_2') }}\right)^2 \leq a\,\left(\frac{|x_2-y_2|^{\phantom{ (\rho,\rho')} }}{
(|x_1|+|y_1|)^{ (\rho,\rho') }}\right)^2\le a
\end{eqnarray*}
for all  $|x_2-y_2|\le |x_1|+|y_1|^{(\rho,\rho')} $.
Finally
\[
 c_{\delta_2,\delta_2'}(x_2)\left(\nabla_{x_2}
\frac{|x_2-y_2|^{\phantom{(\delta_2,\delta_2')}} }{  (|x_1|+|y_1|)^{(\delta_2,\delta_2')}   }\right)^2
\leq \left(
\frac{|x_1|^{(\delta_2,\delta_2')}  }{
(|x_1|+|y_1|)^{(\delta_2,\delta_2')}   }\right)^2
\le a
\;.
\]
This completes the verification of (\ref{csg.nabla}).

  Since $D_\delta$ satisfies (\ref{csg.nabla}) it follows formally from the definition  (\ref{ecsg4.00})
  of the Riemannian distance that $D_\delta (x_1,x_2\,;y_1,y_2) \leq a\, d_\delta (x_1,x_2\,;y_1,y_2)$.
 One cannot, however, immediately make this deduction since $D_\delta\not\in W^{1,\infty}(\Ri^n\times \Ri^m)$. 
  It is, however, a  continuous    function which is locally Lipschitz differentiable on
$(\Ri^n -   \{ 0  \} )\times \Ri^m$.
But the distance is not changed if one replaces the space of trial functions
$ W^{1,\infty}(\Ri^n\times \Ri^m)$ in the definition   (\ref{ecsg4.00}) by a space of functions which are Lipschitz differentiable on the complement  of a closed set of measure zero.
This can be deduced by remarking that since we may assume the coefficients are continuous both 
definitions agree with the shortest path definition of the distance. 
Therefore $D_\delta (x_1,x_2\,;y_1,y_2) \leq a\, d_\delta (x_1,x_2\,;y_1,y_2)$ and since we have already established the converse inequality one concludes that
$D_\delta (x_1,x_2\,;y_1,y_2) \sim d_\delta (x_1,x_2\,;y_1,y_2)$.

\smallskip

It remains to  prove the volume estimates (\ref{evol}).
The proof will be divided into three steps.
First we consider small $r$, small compared with $|x_1|$, secondly we consider 
large $r$ and finally we deal with  intermediate values.

\smallskip

\noindent{\bf Step 1}$\;$ Assume $c\,r \leq |x_1|^{(1-\delta_1,1-\delta_1')}$ with $c>1$.
In fact we will choose $c\gg1$ in the course of the proof.
We now argue that there are $a_1, a_2\in\Ri$ with $0<a_1<a_2$ such that 
\begin{eqnarray}
(x_1,x_2)&&+
\left[-a_1\,r\,|x_1|^{(\delta_1,\delta_1')}, a_1\,r\,|x_1|^{(\delta_1,\delta_1')}\right]^n
\times
\left[-a_1\,r\,|x_1|^{(\delta_2,\delta_2')}, a_1\,r\,|x_1|^{(\delta_2,\delta_2')}\right]^m
\subset  B(x_1,x_2\,;r)\nonumber \\[5pt]
&&\hspace{-1cm}{}\subset
(x_1,x_2)+
\left[-a_2\,r\,|x_1|^{(\delta_1,\delta_1')}, a_2\,r\,|x_1|^{(\delta_1,\delta_1')}\right]^n
\times
\left[-a_2\,r\,|x_1|^{(\delta_2,\delta_2')}, a_2\,r\,|x_1|^{(\delta_2,\delta_2')}\right]^m\;\;.
\label{incl1}
\end{eqnarray}
Once this is established one has the volume estimates
\begin{equation}
|B(x_1,x_2\,;r)|\sim r^{n+m}|x_1|^{(n\delta_1+m\delta_2,n\delta_1'+m\delta_2')}
=r^{n+m}|x_1|^{(\beta,\beta')}
\label{ecsg5.101}
\end{equation}
for $r \leq c^{-1}\, |x_1|^{(1-\delta_1,1-\delta_1')}$

First consider the right hand inclusion of (\ref{incl1}).
Set $\|x\|=\max_{1\leq k\leq n}|x^{(k)}|$ and $\|y\|=\max_{1\leq l\leq m}|y^{(l)}|$
where $x^{(k)}$ and $y^{(l)}$ are the components of $x\in\Ri^n$ and $y\in \Ri^m$, respectively.
Thus we  have to prove that if $(y_1,y_2)\in B(x_1,x_2\,;r)$ then 
$\|x_1-y_1\|\leq a_2\,r\,|x_1|^{(\delta_1,\delta_1')}$ and 
$\|x_2-y_2\|\leq a_2\,r\,|x_1|^{(\delta_2,\delta_2')}$.
But $d_\delta (x_1,x_2;y_1,y_2)\sim D_\delta (x_1,x_2;y_1,y_2)$.
Hence if $(y_1,y_2)\in B(x_1,x_2\,;r)$ then 
$\|x_1-y_1\|\leq |x_1-y_1|\leq a\,r\,(|x_1|+|y_1|)^{(\delta_1,\delta_1')}$.
Therefore
\begin{eqnarray*}
|\|x_1\|-\|y_1\||\leq \|x_1-y_1\|&\leq &a'\,c^{-1}\,(|x_1|+|y_1|)^{(1-\delta_1,1-\delta_1')}\,(|x_1|+|y_1|)^{(\delta_1,\delta_1')}\\[5pt]
&=&a'\,c^{-1}\,(|x_1|+|y_1|)\leq a'\,c^{-1}\,(n+m)(\|x_1\|+\|y_1\|)
\;.
\end{eqnarray*}
Choosing $c$ large one deduces that $\|y_1\|\sim\|x_1\|$ and $|y_1|\sim |x_1|$.
In particular one has 
\[
\|x_1-y_1\|\leq a\,r\,(|x_1|+|y_1|)^{(\delta_1,\delta_1')}\leq a_2\,r\,|x_1|^{(\delta_1,\delta_1')}
\]
as required.

Secondly, since  $(y_1,y_2)\in B(x_1,x_2\,;r)$ one has $ \Delta_\delta (x_1,x_2\,;y_1,y_2)\leq a\,r$.
There are two cases to consider.  
The first is if $ |x_2-y_2|\geq (|x_1|+|y_1|) ^{(\rho,\rho')}$  then   $\Delta_\delta (x_1,x_2\,;y_1,y_2)= 
|x_2-y_2|^{(1-\gamma,1-\gamma')}  $.
Hence 
\[
|x_2-y_2|^{(1-\gamma,1-\gamma')}  \leq a\,r\leq a\,c^{-1}\,|x_1|^{(1-\delta_1,1-\delta_1')}
\leq a\,c^{-1}\,(|x_1|+|y_1|)^{(1-\delta_1,1-\delta_1')}
\;.
\]
Thus if $c\geq a$ then 
\[
|x_2-y_2|\leq\,(|x_1|+|y_1|)^{((1-\delta_1)(1-\gamma)^{-1},(1-\delta_1')(1-\gamma')^{-1})}
=(|x_1|+|y_1|)^{(\rho,\rho')}
\]
which is in contradiction with the assumption  $ |x_2-y_2|\geq (|x_1|+|y_1|) ^{(\rho,\rho')}$.
Therefore one must have  $ |x_2-y_2|\leq (|x_1|+|y_1|) ^{(\rho,\rho')}$ and 
$\Delta_\delta (x_1,x_2\,;y_1,y_2)= {|x_2-y_2|}/{(|x_1|+|y_1|)^{(\delta_2,\delta_2')}}$.
Then, however,
\[
\|x_2-y_2\|\leq |x_2-y_2|\leq a\,r\,(|x_1|+|y_1|)^{(\delta_2,\delta_2')}\leq a_2\,r\,|x_1|^{(\delta_2,\delta_2')}
\]
because $|y_1|\sim|x_1|$ by the previous argument.
Thus the proof of the right hand inclusion of (\ref{incl1}) is complete.

Thirdly, consider the left hand inclusion of (\ref{incl1}).
Now we need to prove that if one has 
$\|x_1-y_1\|\leq a_1\,r\,|x_1|^{(\delta_1,\delta_1')}$ and 
$\|x_2-y_2\|\leq a_1\,r\,|x_1|^{(\delta_2,\delta_2')}$ then 
$(y_1,y_2)\in B(x_1,x_2\,;r)$.
But by the first assumption
\[
|x_1-y_1|/(|x_1|+|y_1|)^{(\delta_1,\delta_1')}\leq  a_1\,n\,r\,|x_1|^{(\delta_1,\delta_1')}/(|x_1|+|y_1|)^{(\delta_1,\delta_1')}\leq a_1\,n\,r
\;.
\]
Then by the second assumption
\begin{eqnarray*}
|x_2-y_2|\leq a_1\,m\,r\,|x_1|^{(\delta_2,\delta_2')}
&\leq& a_1\,m\,c^{-1}\,|x_1|^{(1+\delta_2-\delta_1,1+\delta_2'-\delta_1')}\\[5pt]
&=&a_1\,m\,c^{-1}\,|x_1|^{(\rho,\rho')}
\leq a_1\,m\,c^{-1}\,(|x_1|+|y_1|)^{(\rho,\rho')}
\;.
\end{eqnarray*}
Hence if $c\geq a_1\,m$ then $|x_2-y_2|\leq  (|x_1|+|y_1|)^{(\rho,\rho')}$.
Therefore it follows that  $\Delta_\delta (x_1,x_2\,;y_1,y_2)= {|x_2-y_2|}/{(|x_1|+|y_1|)^{(\delta_2,\delta_2')}}$.
Hence using $\|x_2-y_2\|\leq a_1\,r\,|x_1|^{(\delta_2,\delta_2')}$  again one has
\[
\Delta_\delta (x_1,x_2\,;y_1,y_2)\leq  a_1\,m\,r\,|x_1|^{(\delta_2,\delta_2')}/{(|x_1|+|y_1|)^{(\delta_2,\delta_2')}}\leq a_1\,m \,r
\]
and 
\[
D_\delta (x_1,x_2\,;y_1,y_2)=
{|x_1-y_1|}/{(|x_1|+|y_1|)^{(\delta_1,\delta_1')}}
+\Delta_\delta (x_1,x_2\,;y_1,y_2)\leq a_1\,(n+m)\,r
\;.
\]
Hence if $a_1$ is sufficiently small one concludes that $(y_1,y_2)\in B(x_1,x_2\,;r)$.

\smallskip

\noindent{\bf Step 2}$\;$ Assume $r/c \geq |x_1|^{(1-\delta_1,1-\delta_1')}$
with $c>1$ where we will again choose $c\gg1$.
We now argue that there are
$a_1, a_2\in\Ri$ with $0<a_1<a_2$ such that 
\begin{eqnarray}
(0,x_2)+
\left[-a_1\,r^{(\sigma,\sigma')},
a_1\,r ^{(\sigma,\sigma')}\right]^n
&\times&
\left[-a_1\,r ^{(\rho\sigma,\rho'\sigma')},
a_1\,r ^{(\rho\sigma,\rho'\sigma')} \right]^m
\subset B(x_1,x_2\,;r)\nonumber\\[5pt]
&&\hspace{-4.5cm}{}\subset
(0,x_2)+
\left[-a_2\,r^{(\sigma,\sigma')},
a_2\,r ^{(\sigma,\sigma')}\right]^n
\times
\left[-a_2\,r ^{(\rho\sigma,\rho'\sigma')},
a_2\,r ^{(\rho\sigma,\rho'\sigma')}\right]^m
\label{incl2}
\end{eqnarray}
where $\sigma=(1-\delta_1)^{-1}$ and $\sigma'=(1-\delta'_1)^{-1}$.
These inclusions then yield the volume estimates
\begin{equation}
|B(x_1,x_2\,;r)|\sim r^{(n\sigma, n\sigma')}r^{(m\rho\sigma, m\rho'\sigma')}=r^{(D,D')}
\;.
\label{ecsg5.102}
\end{equation}

The first step in deducing  (\ref{incl2}) is to observe that if $c$ is sufficiently large
then
\[
B(0,x_2\,;r/2)\subseteq B(x_1,x_2\,;r)\subseteq B(x_1,x_2\,;2r)
\]
and so the proof is effectively reduced to the case $x_1=0$.
Then the rest of the proof is similar to the argument  in Step~1 but somewhat simpler
because of the choice $x_1=0$.

\smallskip

\noindent{\bf Step 3}$\;$
It follows from Step~1, and in particular (\ref{ecsg5.101}), that the volume estimates (\ref{evol}) are valid
if $ |x_1|^{(1-\delta_1,1-\delta_1')}\geq c\, r$.
Alternatively, it follows from Step~2, and in particular (\ref{ecsg5.102}),
that the estimates  are valid if  $r \geq c\, |x_1|^{(1-\delta_1,1-\delta_1')}$.
Therefore we now assume that $r\in[c^{-1} |x_1|^{(1-\delta_1,1-\delta_1')},c \,|x_1|^{(1-\delta_1,1-\delta_1')}]$.
But then setting $r_1=r/c$ and $r_2=r\,c$ one has
 $B(x_1,x_2\,;r_1) \subset B(x_1,x_2\,;r) \subset B(x_1,x_2\,;r_2)$.
 But $r_1\leq  |x_1|^{(1-\delta_1,1-\delta_1')}$.
 Hence
 $B(x_1,x_2\,;r_1)~\sim r_1^{n+m}|x_1|^{(\beta,\beta')}\sim
 r^{n+m}\,|x_1|^{(\beta,\beta')}$ by Step~1.
 Moreover, $r_2\geq|x_1|^{(1-\delta_1,1-\delta_1')}$.
 Hence 
$B(x_1,x_2\,;r_2)~\sim r_2^{(D,D')} \sim
 r^{(D,D')}$
 by Step~2.
 Combining these estimates one concludes that 
 \[
 a\,r^{n+m}\,|x_1|^{(\beta,\beta')}\leq |B(x_1,x_2\,;r) |\leq a'\,r^{(D,D')}
 \]
 for all $r\in[c^{-1} |x_1|^{(1-\delta_1,1-\delta_1')},c \,|x_1|^{(1-\delta_1,1-\delta_1')}]$.
But in this range $r^{(D,D')}\sim a\,r^{n+m}\,|x_1|^{(\beta,\beta')}$.
Therefore  the volume estimates are established.
 \hfill$\Box$
 
 \ruimte
 
 The volume estimates allow one to prove the doubling property and to identify the doubling dimension.
 
 \begin{cor}\label{ccsg5.1} The Riemannian balls $B(x_1,x_2\,;r)$ associated with the Gru\v{s}in operator satisfy the doubling property
 \[
 |B(x_1,x_2\,;s\,r)|\leq a\,s^{(D\vee D')} |B(x_1,x_2\,;r)|
 \]
 for all $(x_1,x_2)\in \Ri^{n+m}$ and all $s\geq 1$.
 \end{cor}
 \proof\
 There are three cases to consider.
 
 \smallskip
 
\noindent{\bf Case~1}$\;\;r\leq r\,s\leq |x_1|^{(1-\delta_1,1-\delta_1')}$.
Then the volume estimates of Proposition~\ref{pcsg5.1} give
 \[
 |B(x_1,x_2\,;s\,r)|\sim(r\,s)^{n+m}\,|x_1|^{(\beta,\beta')}\sim s^{n+m}\, |B(x_1,x_2\,;r)|
 \;.
 \]
 But $n+m\leq D\vee D'$ and $s\geq 1$ so $s^{n+m}\leq s^{(D\vee D')}$ and the doubling property follows.

\smallskip

\noindent{\bf Case~2}$\;\;|x_1|^{(1-\delta_1,1-\delta_1')}\leq r\leq rs$.
Then the volume estimates  give
 \[
 |B(x_1,x_2\,;s\,r)|\sim(r\,s)^{(D,D')}\leq s^{(D\vee D')}\, r^{(D,D')}\sim s^{(D\vee D')}\,  |B(x_1,x_2\,;r)|
\]
  and the doubling property is established.
  
\smallskip

\noindent{\bf Case~3}$\;\;r\leq |x_1|^{(1-\delta_1,1-\delta_1')}\leq  rs$.
Then the volume estimates  give
 \[
 |B(x_1,x_2\,;s\,r)|\sim(r\,s)^{(D,D')}\leq s^{(D\vee D')}\, r^{(D,D')}
 \;.
\]
But
\[
r^{(D,D')}=r^{n+m}r^{(\beta(1-\delta_1)^{-1},\beta'(1-\delta_1')^{-1})}\leq 
r^{n+m}|x_1|^{(\beta,\beta')}\sim |B(x_1,x_2\,;r)|
\;.
\]
Combination of these estimates gives the doubling property again.
 \hfill$\Box$

\bigskip

The doubling property has been established for a different class of Gru\v{s}in operators in
\cite{FrS} (see also \cite{Fra}  \cite{FGW2}).

\section{Kernel bounds}\label{Scsg6} 

In this section we prove that the  semigroup kernel $K_t$ associated with a
general Gru\v{s}in  operator $H$
 conserves probability, i.e.\  the evolution is stochastically complete.
Then  we establish that $K_t$  satisfies off-diagonal volume dependent upper bounds.
Subsequently it is  possible to apply standard reasoning to obtain further more detailed properties 
such as Gaussian bounds and on-diagonal lower bounds.
This will be discussed at the end of the section.

First, the semigroup $S$ generated by the Gru\v{s}in operator $H$ is submarkovian. 
In particular it extends to a contractive 
semigroup on $L_\infty(\Ri^n\times\Ri^m)$.
Therefore 
\[
0\leq \esssup_{x\in\Ri^{n+m}}\int_{\Ri^{n+m}}dy\,K_t(x\,;y)\leq 1
\;.
\]
But the particular structure of the operator gives a stronger result.

\begin{thm}\label{tcsg6.1} 
The semigroup $S$ associated with the Gru\v{s}in operator  on $L_\infty(\Ri^n\times\Ri^m)$ satisfies
$S_t\one=\one$ for all $t>0$.
Hence the kernel satisfies
\begin{equation}
\int_{\Ri^n\times\Ri^m}dy\,K_t(x\,;y)=1
\label{ecsg6.1}
\end{equation}
for  all $t>0$ almost all $x\in \Ri^{n+m}$.
\end{thm}
\proof\
Given that the semigroup satisfies the $L_2$ off-diagonal bounds of Proposition~\ref{pscsg4.0}
and the volume of the Riemannian balls have polynomial growth, by the estimates of 
Proposition~\ref{pcsg5.1}, one can prove the theorem by  a slight variation of the argument given in Proposition~3.6 of \cite{ERSZ1} but estimating with respect to a Riemannian distance instead of the Euclidean distance.
We omit the details.
\hfill$\Box$

\bigskip

An alternative proof of the theorem can be constructed by  approximating $H$ through strongly elliptic operators.
Let  $H_{N,\varepsilon}$ be  the strongly
elliptic approximants,  with coefficients $C_{N,\varepsilon}=(C\wedge NI)+\varepsilon I$, to the Gru\v{s}in operator.
Then the  semigroups $S^{(N,\varepsilon)}_t$ generated by the   $H_{N,\varepsilon}$
converge strongly on $L_2$ to the semigroup $S_t$ in the double limit $N\to\infty$ followed by $\varepsilon\to 0$.
But the convergence is stronger.

\begin{prop}\label{pcsg4.9}
The semigroups $S^{(N,\varepsilon)}_t$ 
converge strongly to  $S_t$ on each of the  $L_p$-spaces with $p\in[1,\infty\rangle$
and in the weak$^*$ sense on $L_\infty$.
\end{prop}

Since this result is not used in the sequel we omit the proof.
But once it is established then Theorem~\ref{tcsg6.1} follows since
 $S^{(N,\varepsilon)}_t\one=\one$ because the  $H_{N,\varepsilon}$ are strongly elliptic.
Hence $S_t\one=\one$ in the  weak$^*$ limit.

\begin{remarkn}\label{rcsg6.10}
The proof of Proposition~\ref{pcsg4.9} and Theorem~\ref{tcsg6.1} 
 uses very little structure of the Gru\v{s}in operator.
 It only requires $L_2$ off-diagonal bounds and polynomial volume growth.
 The first are given for general elliptic operators by Proposition~\ref{pscsg4.0} and the polynomial growth
 follows from the Gru\v{s}in structure.
 It follows from this observation that similar statements would also follow for the degenerate operators considered in \cite{FrS} of \cite{FGW2}.
 \end{remarkn}

Next we consider upper bounds on the semigroup kernel.

\begin{thm}\label{tcsg6.2} 
There is an $a>0$ such that the  semigroup kernel $K$ of the Gru\v{s}in operator $ H$
satisfies 
 \begin{equation}
0\leq  K_t(x\,;y)\leq a\,(|B(x\,;t^{1/2})|\,|B(y\,;t^{1/2})|)^{-1/2}
\label{ecsg6.2}
\end{equation}
for  all $t>0$ and  almost all $x,y\in \Ri^{n+m}$.
\end{thm}

Since the semigroup $S$ is self-adjoint it follows that the semigroup kernel is positive-definite.
Therefore one formally has $|K_t(x\,;y)|^2\leq K_t(x\,;x)K_t(y\,;y)$ and the estimate
apparently reduces to an on-diagonal estimate.
But this is only a formal calculation since the kernel is not necessarily continuous and its diagonal value
is not necessarily defined.
Nevertheless the starting point of the proof is a set-theoretic reduction to an on-diagonal estimate.

\begin{lemma}\label{lcsg6.1}
Let $X,Y$ be open sets and define $K_t(X\,;Y)=\esssup_{x\in X,y\in Y}K_t(x\,;y)$.
Then
\[
|K_t(X\,;Y)|^2\leq K_t(X\,;X)\,K_t(Y\,;Y)
\;.
\]
\end{lemma}
\proof\
First observe that 
\begin{eqnarray*}
K_t(X\,;Y)=\|\one_YS_t\one_X\|_{1\to\infty}
\leq\|S_{t/2}\one_X\|_{1\to2}\|\one_Y S_{t/2}\|_{2\to\infty}
= \|\one_XS_{t/2}\|_{2\to\infty}\|\one_Y S_{t/2}\|_{2\to\infty}
\;.
\end{eqnarray*}
But if $T$ is  bounded from $L_2$ to $L_\infty$ then $(\|T\|_{2\to\infty})^2=\|TT^*\|_{1\to\infty}$.
Therefore
\begin{eqnarray*}
|K_t(X\,;Y)|^2
\leq\|\one_XS_{t}\one_X\|_{1\to\infty}\|\one_Y S_{t}\one_Y\|_{1\to\infty}
=K_t(X\,;X)\,K_t(Y\,;Y)
\end{eqnarray*}
as required.\hfill$\Box$

\bigskip

\noindent{\bf Proof of Theorem~\ref{tcsg6.2} } It follows from the lemma that 
\[
K_t(x\,;y)\leq \inf_{X\ni x}K_t(X\,;X)^{1/2}\,\inf_{Y\ni y}K_t(Y\,;Y)^{1/2}
\;.
\]
Thus it suffices to prove that 
\[
 \inf_{X\ni x}K_t(X\,;X)\leq a\,|B(x\,;t^{1/2})|^{-1}
 \;.
 \]
 There are two distinct cases corresponding to the different volume behaviours given by 
 Proposition~\ref{pcsg5.1}.
 
 First, let  $x=(x_1,x_2)$ and suppose $|x_1|^{(1-\delta_1,1-\delta_1')}\leq t^{1/2}$. 
 Then 
 $|B(x\,;t^{1/2})|\sim t^{(D/2,D'/2)}$, by  Proposition~\ref{pcsg5.1},
and $\|S_t\|_{1\to\infty}\leq a\,t^{(-D/2,-D'/2)}$, by Proposition~\ref{pcsg3.0}.
 Therefore 
 \[
 \inf_{X\ni x}K_t(X\,;X)\leq \|K_t\|_\infty=\|S_t\|_{1\to\infty}\leq a\,t^{(-D/2,-D'/2)}
 \leq a'\,|B(x\,;t^{1/2})|^{-1}
 \;.
 \]
 
Secondly, suppose that $t^{1/2}\leq |x_1|^{(1-\delta_1, 1-\delta_1')}$. 
 Then  $|B(x\,;t^{1/2})|\sim  t^{(n+m)/2}|x_1|^{(\beta,\beta')}$ by  Proposition~\ref{pcsg5.1}.
 This case is considerably more difficult to analyze and it is here that we apply the comparison techniques of Section~\ref{Scsg4}.
 
Set $r=|x_1|$.
Let $C$ denote the coefficient matrix of $H$ and  choose $a_n, a_m>0$ such that 
$C(y)\geq a_n\,|y_1|^{(2\delta_1,2\delta_1')} I_n+a_m\,|y_1|^{(2\delta_2,2\delta_2')} I_m$ for all $y=(y_1,y_2)$ with $|y_1|\leq r/2$.
Next set $C_r(y)=C(y)$ if $|y_1|>r/2$ and
$C_r(y)= a_n\,r^{(2\delta_1,2\delta_1')} I_n+a_m\,r^{(2\delta_2,2\delta_2')} I_m$ if $|y_1|\leq r/2$.
Then $C_r\geq a\,C$ for a suitable $a>0$. 
Let  $H_r$ be  the Gru\v{s}in operator with coefficient matrix $C_r$ 
and set $H_1=H_r$ and $H_2=a\,H$.
Then $H_1\geq H_2$.
Moreover $H_1\geq \mu I>0$  for some 
$\mu=a\,(a_nr^{(2\delta_1,2\delta_1')}\wedge a_m r^{(2\delta_2,2\delta_2')})$.
Thus the basic assumptions of Theorem~\ref{tcsg3.1} and Corollary~\ref{ccsg3.1} are satisfied with
this choice.
Now $H_r$ is a  Gru\v{s}in operator  with local parameters $\delta_1=0=\delta_2$ but with the same
global parameters $\delta_1',\delta_2'$ as $H$.
Therefore
\[
\|S^{(1,0)}_t\|_{1\to\infty}\vee \|S^{(2,0)}_t\|_{1\to\infty}\leq V(t)^{-1}
\]
for all $t>0$ where $V(t)=a\,t^{(D/2,D'/2)}$.
Since $V$ satisfies the doubling property one may now apply Corollary~\ref{ccsg3.1}.
Note that for this application $U=\{y: |y_1|\leq r/2\}$ and $d(x\,;U)\sim r^{(1-\delta_1,1-\delta_1')}$.

Next we need an improved estimate on the crossnorm $\|S^{(1,0)}_t\|_{1\to\infty}$ of the comparison 
semigroup.
Let $\widehat H$ denote the constant coefficient operator with coefficients $\widehat C_r
= a_n\,r^{(2\delta_1,2\delta_1')} I_n+a_m\,r^{(2\delta_2,2\delta_2')} I_m$.
Then $\widehat H$ is a Fourier multiplier and the corresponding function $F$  is given by
 $F(p_1,p_2)=a_n\,r^{(2\delta_1,2\delta_1')}p_1^2+a_m\,r^{(2\delta_2,2\delta_2')} p_2^2$.
But there is an $a>0$ such that $H_1 \geq a\, \widehat H$ and we can apply Lemma~\ref{lcsg2.2}
to obtain a uniform bound on $ \|S^{(1,0)}_t\|_{1\to\infty}$.
First, suppose $r=1$ then one immediately has $ \|S^{(1,0)}_t\|_{1\to\infty}\leq a\,t^{-(n+m)/2}$
for all $t>0$.
Secondly, the introduction of $r$ corresponds to a dilation of $\Ri^{n}\times \Ri^{m}$
with each direction in $\Ri^{n}$ dilated by $r^{(-\delta_1,-\delta_1')}$  and each direction in 
$\Ri^{m}$ dilated by $r^{(-\delta_2,-\delta_2')}$.
The dilation adds a factor to the crossnorm corresponding to the Jacobian
$r^{(-n\delta_1-m\delta_2,-n\delta_1'-m\delta_2')}=r^{(-\beta,-\beta')}$ of the dilation. 
Therefore  
\[
 \|S^{(1,0)}_t\|_{1\to\infty}\leq a'\,t^{-(n+m)/2}r^{(-\beta,-\beta')}
 \leq a''\,|B(x\,;t^{1/2})|^{-1}
 \]
 where the last bound follows from the second estimate of Proposition~\ref{pcsg5.1}.

Now we may apply Corollary~\ref{ccsg3.1} with $A=X$.
One obtains an estimate
 \begin{eqnarray*}
 \inf_{X\ni x}K_t(X\,;X)= \inf_{X\ni x}K^{(2,0)}_t(X\,;X)
 &\leq&
 \|S^{(1,0)}_t\|_{1\to\infty}+
\mathop{\smash \inf\vphantom\sup}_{X\ni x}\sup_{y,z\in X}
 |K^{(1,0)}_t(y\,;z)-K^{(2,0)}_t(y\,;z)|\\[5pt]
 &\leq& a\,|B(x\,;t^{1/2})|^{-1}+a\,V(t^2/\rho^2)^{-1}\,(\rho^2/t)^{-1/2}
 \,e^{-\rho^2/(4t)}\\[5pt]
 &\leq& a\,|B(x\,;t^{1/2})|^{-1}\left(1+R(x\,;t)\right)
\end{eqnarray*}
where 
\[
R(x\,;t)=|B(x\,;t^{1/2})|\,V(t^2/\rho^2)^{-1}\,(\rho^2/t)^{-1/2}
 \,e^{-\rho^2/(4t)}
 \;,
 \]
 $\rho=d(x\,;U)\sim r^{(1-\delta_1,1-\delta_1')}$ and $V(t)=a\,t^{(D/2,D'/2)}$.
Thus it suffices to show that $R(x\,;t)$ is uniformly bounded for $x,t$ satisfying
 $t^{1/2}\leq r^{(1-\delta_1, 1-\delta_1')}$ with $r=|x_1|$. 
It is necessary to distinguish between three cases.

\smallskip

\noindent{\bf Case~1}$\;\;t,r\leq 1$.
In this case  $t^{1/2}\leq r^{1-\delta_1}\leq 1$ and $\rho\sim r^{1-\delta_1}$. 
Hence $\rho^2/t\sim (r^{1-\delta_1}t^{-1/2})^2$ and  $t^2/\rho^2\leq a\,t$. 
Therefore, since  $t\leq1$, one must have
$V(t^2/\rho^2)^{-1}\leq a\, t^{-D/2}$.
Moreover, $|B(x\,;t^{1/2})|\sim  t^{(n+m)/2}r^{\beta}$.
Then one  computes that 
\begin{eqnarray*}
R(x\,;t)&\leq&
  a\,t^{(n+m)/2}r^{\beta}\,t^{-D/2}\,(r^{1-\delta_1}t^{-1/2})^{-1}e^{-a'(r^{1-\delta_1}t^{-1/2})^2}\\[5pt]
  &=& a\,(r^{1-\delta_1}t^{-1/2})^{-1+\beta(1-\delta_1)^{-1}}\,e^{-a'(r^{1-\delta_1}t^{-1/2})^2}
  \;.
  \end{eqnarray*}
  But since $r^{1-\delta_1}t^{-1/2}\geq1$  it follows that $R(x\,;t)$ is uniformly bounded.

\smallskip

\noindent{\bf Case~2}$\;\;t,r\geq 1$.
In this case  $1\leq t^{1/2}\leq r^{1-\delta_1'}$ and $\rho\sim r^{1-\delta_1'}$. 
Hence $\rho^2/t\sim (r^{1-\delta_1'}t^{-1/2})^2$ and  $t^2/\rho^2\leq a\,t$. 
But now  $t\geq1$ and so 
$V(t^2/\rho^2)^{-1}\leq a\, t^{-D'/2}$.
Moreover, $|B(x\,;t^{1/2})|\sim  t^{(n+m)/2}r^{\beta'}$.
Then one  computes  as in Case~1 that
\begin{eqnarray*}
R(x\,;t)&\leq&
  a\,t^{(n+m)/2}r^{\beta'}\,t^{-D'/2}\,(r^{1-\delta_1'}t^{-1/2})^{-1}e^{-a'(r^{1-\delta_1'}t^{-1/2})^2}\\[5pt]
  &=& 
 a\,(r^{1-\delta_1'}t^{-1/2})^{-1+\beta'(1-\delta_1')^{-1}}\,e^{-a'(r^{1-\delta_1'}t^{-1/2})^2}
  \;.
  \end{eqnarray*}
  But since $r^{1-\delta_1'}t^{-1/2}\geq1$  it again  follows that $R(x\,;t)$ is uniformly bounded.

\smallskip

\noindent{\bf Case~3}$\;\;t\leq 1, r\geq 1$.
This is a hybrid case which is rather different to the previous two cases.
One again has $\rho\sim r^{1-\delta_1'}$ and  $\rho^2/t\sim (r^{1-\delta_1'}t^{-1/2})^2$.
Moreover, $t^2/\rho^2\leq a\,t$. 
But since $t\leq 1$ one has $V(t^2/\rho^2)^{-1}\leq a\, t^{-D/2}$.
In addition $|B(x\,;t^{1/2})|\sim  t^{(n+m)/2}r^{\beta'}$.
Therefore one now estimates that
\begin{eqnarray*}
R(x\,;t)&\leq&
  a\,t^{(n+m)/2}r^{\beta'}\,t^{-D/2}\,(r^{1-\delta_1'}t^{-1/2})^{-1}e^{-a'(r^{1-\delta_1'}t^{-1/2})^2}\\[5pt]
   &= &a\,r^{\beta'}\,(t^{-1/2})^{\beta(1-\delta_1)^{-1}}\,
  (r^{1-\delta_1'}t^{-1/2})^{-1} e^{-a'(r^{1-\delta_1'}t^{-1/2})^2}
  \;.
  \end{eqnarray*}
  But now $r^{1-\delta_1'}t^{-1/2}\geq 1$.
 Consequently for each $N\geq 1$ there is an $a_N$ such that 
\[
R(x\,;t)\leq 
a_N\,r^{\beta'}\,(t^{-1/2})^{\beta(1-\delta_1)^{-1}}\,(r^{1-\delta_1'}t^{-1/2})^{-N}
 \;.
 \]
 Choosing $N$ large ensures that this expression is uniformly bounded
 for all $r\geq1$ and $t\leq 1$.
 
 \smallskip
 
 The proof of the theorem is now complete.\hfill$\Box$

 \bigskip
 
 Theorems~\ref{tcsg6.1}  and \ref{tcsg6.2} have a number of standard implications.
 First one may convert the volume bounds of the latter theorem into  Gaussian bounds.
 
 \begin{cor}\label{ccsg6.1}
For each $\varepsilon>0$ there is an $a>0$ such that the  semigroup kernel $K$ of the Gru\v{s}in operator satisfies
 \begin{equation}
0\leq  K_t(x\,;y)\leq a\,(|B(x\,;t^{1/2})|\,|B(y\,;t^{1/2})|)^{-1/2}\,e^{-d(x;y)^2/(4(1+\varepsilon)t)}
\label{ecsg6.21}
\end{equation}
for  all $t>0$ and almost all $x,y\in \Ri^{n+m}$.
\end{cor}

There are several different arguments for passing from on-diagonal kernel bounds to Gaussian bounds
(see, for example, the lecture notes of Grigor'yan \cite{Gri3}).
One proof of the corollary which is in the spirit of the present paper is given by Theorem~4 of 
\cite{Sik3}.
Note that in the latter reference it is implicitly assumed that the kernel is well-defined on the  diagonal
but this is not essential.
One can argue with open sets and near diagonal estimates as in the proofs of Theorem~\ref{tcsg3.1} and  Theorem~\ref{tcsg6.2}.

\smallskip

It is also a standard argument to pass from the Gaussian bounds of Corollary~\ref{ccsg6.1} and the
conservation property of Theorem~\ref{tcsg6.1} to on-diagonal lower bounds.
Again one has to avoid problems with the definition of the diagonal values.

 \begin{cor}\label{ccsg6.2}
There is an $a>0$ such that the  semigroup kernel $K$ of the Gru\v{s}in operator satisfies
 \begin{equation}
\inf_{X\ni x}|X|^{-2}\int_Xdx\int_Xdy\, K_t(x\,;y)
 \geq a\,|B(x\,;t^{1/2})|^{-1}
 \label{ecsg6.22}
\end{equation}
for  all $t>0$ and almost all $x\in \Ri^{n+m}$ where the average is over open subsets $X$.
\end{cor}
\proof\
At the risk of confusion with the earlier definition we set 
\[
K_t(X\,;Y)=\int_Xdx\int_Ydy\, K_t(x\,;y)=(\one_X,S_t\one_Y)
\]
for each pair of bounded open sets $X,Y$.
Then using self-adjointness and the semigroup property one again verifies that 
\[
|K_t(X\,;Y)|^2\leq K_t(X\,;X)\,K_t(Y\,;Y)
\;.
\]
Now fix $x$ and $X\ni x$ and  let $Y=B(x\,;R\,t^{1/2})$ for  $R>0$.
Then 
\[
K_t(X\,;Y)=\int_Xdx\Big(1-\int_{Y^{\rm c}}dy\, K_t(x\,;y)\Big)
\]
for all $t>0$ by Theorem~\ref{tcsg6.1}.
But now using the Gaussian bounds of Theorem~\ref{tcsg6.2} and choosing $R$ large one can ensure that
\[
K_t(X\,;Y)\geq |X|/2
\;.
\]
Then, however, with this choice of $Y$ one has 
\[
|X|^{-2}K_t(X\,;X)\geq (4K_t(Y\,;Y))^{-1}
\;.
\]
But using the bounds of Theorem~\ref{tcsg6.2} one immediately finds 
\[
K_t(Y\,;Y)\leq a\,\Big(\int_{d(x;y)<Rt^{1/2}}dy\,|B(y\,;t^{1/2})|^{-1/2}\Big)^2 
\;.
\]
If $y\in B(x\,;R\,t^{1/2})$ the doubling property gives
 \[
 |B(x\,;t^{1/2})|\leq  |B(y\,;d(x\,;y)+t^{1/2})|\leq  |B(y\,;(1+R)t^{1/2})|
 \leq a\,(1+R)^{\tilde D} |B(y\,;t^{1/2})|
 \]
with $\tilde D=D\vee D'$.
 In addition 
 $|B(x\,;R\,t^{1/2})|\leq a\,R^{\tilde D}  |B(x\,;t^{1/2})|.$
 Therefore
 \begin{eqnarray*}
K_t(Y\,;Y)&\leq &a\,(1+R)^{2{\tilde D} }\Big(\int_{d(x;y)<Rt^{1/2}}dy\,|B(x\,;t^{1/2})|^{-1/2}\Big)^2 \\[5pt]
&\leq& a\,(1+R)^{2{\tilde D} }|B(x\,;t^{1/2})|^{-1} |B(x\,;Rt^{1/2})|^2
\leq a'\,((1+R)R)^{2{\tilde D} }|B(x\,;t^{1/2})|
\;.
\end{eqnarray*}
 Combining these estimates gives 
 \[
 |X|^{-2}\int_Xdx\int_Xdy\, K_t(x\,;y)
 \geq a\,|B(x\,;t^{1/2})|^{-1}
 \]
 for all bounded open sets $X$ containing $x$ and this gives the statement of the corollary.
 \hfill$\Box$

\begin{remarkn} \label{rcsg6.1}
If the kernel is continuous  the on-diagonal values are well defined.
Therefore
 \[
K_t(x\,;x)=\inf_{X\ni x}|X|^{-2}\int_Xdx\int_Xdy\, K_t(x\,;y)
 \geq a\,|B(x\,;t^{1/2})|^{-1}
\]
for all $x\in \Ri^{n+m}$ and $t>0$.
\end{remarkn}

\begin{remarkn} \label{rcsg6.100}
The semigroup kernel $K_t$ of the Gru\v{s}in operator $H$ is not necessarily continuous.
In particular if  $n=1$ and $\delta_1\in[1/2,1\rangle$ then the kernel is  discontinuous. 
The discontinuity is a direct consequence of the fact established in 
\cite{ERSZ1} that the action of the corresponding semigroup $S_t$ on $L_2(\Ri\times \Ri^m)$ is not ergodic.
If $H_+=\{(x_1,x_2)\in\Ri\times \Ri^m: x_1>0\}$ and $H_-=\{(x_1,x_2)\in\Ri\times \Ri^m: x_1<0\}$ then
$S_tL_2(H_+)\subseteq L_2(H_+)$ and $S_tL_2(H_-)\subseteq L_2(H_-)$ for all $t>0$.
This is established as follows.

Let $\varphi\in D(h)\subseteq W^{1,2}(\Ri)$ and  set $\varphi_n=\chi_n\varphi$ with  $\chi_n \colon \Ri \to [0,1]$  defined by
\begin{equation}
\chi_n(x)
= \left\{ \begin{array}{ll}
  0 & \mbox{if } x\leq -1 \; ,  \\
  {-\log |x|}/{\log n} & \mbox{if } x\in\langle-1, -n^{-1}\rangle \; , \\
  1 & \mbox{if } x \geq -n^{-1} \;.\label{ecdeg2.137}
         \end{array} \right.
\end{equation}
Then one verifies that $\|\varphi_n-\one_+\varphi\|_2\to 0$ as $n\to\infty$, where $\one_+$ is the indicator function
of $\overline H_+$, and
$h(\varphi_n-\varphi_m)\to0$ as $n,m\to\infty$.
Thus $\one_+\varphi\in D(\overline h)$.
(See \cite{ERSZ1}, Proposition~6.5 and the discussion in  Section~4 of \cite{RSi}.
Note that it is crucial for the last limit that $\delta\in[1/2,1\rangle$.
The conclusion is not valid for 
$\delta\in[0,1/2\rangle$.)
Similarly, replacing $\chi_n(x)$ by $\chi_n(-x)$ one can conclude that 
$\one_{-}\varphi\in D(\overline h)$ with $\one_-$  the indicator function
of $\overline H_-$.
Moreover, $\overline h(\varphi)=\overline h(\one_+\varphi)+\overline h(\one_-\varphi)$.
This suffices to deduce that the semigroup leaves the subspaces $L_2(H_\pm)$ invariant
(see \cite{FOT}, Theorem~1.6.1).

Now one can deduce by contradiction that the kernel has a discontinuity. 
Suppose the kernel $K_t$ of $S_t$ is continuous.
Then it follows from Remark~\ref{rcsg6.1}
that $K_t(x\,;x)\geq a_t>0$ where $a_t=a\,(\sup_{x\in \Ri\times\Ri^m}|B(x\,;t^{1/2}|)^{-1}$.
But since $(\varphi,S_t\psi)=0$ for $\varphi\in L_2(H_+)$ and $\psi\in L_2(H_-)$ one must have
$K_t(x\,;y)=0$ for all $x\in H_+$ and $y\in H_-$.
But this contradicts the continuity hypothesis.
\end{remarkn}

The separation phenomenon in this example  raises the question of boundary conditions  on the hypersurface of separation $x_1=0$.
The closed form $\overline h$ automatically has the decomposition $\overline h(\varphi)=\overline h(\one_+\varphi)+\overline h(\one_-\varphi)$ for all $\varphi\in D(h)$, and by 
closure for all $\varphi\in D(\overline h)$.
Since this is the direct analogue  of the decomposition of the form corresponding to the Laplacian
with Neumann boundary conditions on the hypersurface $x_1=0$ it is tempting to describe the separation in terms of Neumann boundary 
conditions.
But this decomposition is misleading since Neumann and  Dirichlet boundary conditions
coincide in this case.
This can be established by potential theoretic reasoning \cite{RSi}.

  First introduce the  form $h_{\rm Dir}$ by restriction of $h$ to the subspace 
  $D(h)\cap L_{2,c}(H_+\oplus H_-)$ where $L_{2,c}(\Omega)$ is defined as the subspace of $L_2(\Omega)$ spanned
  by the functions with compact support.
  Then  $h_{\rm Dir}$ is closable and its closure $\overline h_{\rm Dir}$ corresponds to the operator $H$ with Dirichlet
  boundary conditions imposed at the boundary $x_1=0$.
 Moreover,  $\overline h_{\rm Dir}\geq \overline h$ in the sense of the ordering of forms.
  But comparison of $\overline h_{\rm Dir}$ and $\overline h$ gives a sharp distinction between the
  weakly  and strongly degenerate cases.  
  
 \begin{prop}\label{psep}
 Consider the Gru\v{s}in operator with $n=1$.
 If $\delta_1\in[0,1/2\rangle$ then  $\overline h_{\rm Dir}> \overline h$ but if  $\delta_1\in[1/2,1\rangle$ then 
 $\overline h_{\rm Dir}=\overline h$. 
\end{prop}
  \proof\
  First, observe that $D(h)\cap L_{2,c}(H_+\oplus H_-)$ is a core of 
  $\overline h_{\rm Dir}$, by definition, and $\overline h_{\rm Dir}=\overline h$ in restriction to the core.
Secondly, it follows from \cite{RSi}, Proposition~3.2, that  $D(h)\cap L_{2,c}(H_+\oplus H_-)$ is a core of 
  $\overline h$ if, and only if, $C_{\overline h}(\{x_1=0\})=0$ where $C_{\overline h}(A)$ denotes the capacity, with
  respect to $\overline h$ of the measurable set $A$.
  Thus $\overline h_{\rm Dir}=\overline h$ if, and only if, $C_{\overline h}(\{x_1=0\})=0$.
  
Now suppose $\delta\in[1/2,1\rangle$ and define
$\xi_n$  by $\xi_n(x_1)=\chi_n(x_1)\wedge \chi_n(-x_1)$
where $\chi_n$ is given by (\ref{ecdeg2.137}).
Then one verifies that $\xi_n\in [0,1]$ and  $\xi_n(x)=1$  for $x\in[-n^{-1},n^{-1}]$.
Moreover, if $\varphi\in D(h)$ then $\varphi_n=\xi_n\varphi\in D(h)$ and 
 $h(\varphi_n)+\|\varphi_n\|_2^2\to0$
as $n\to\infty$.
But this means that $C_{\overline h}(A)=0$  for each bounded measurable subset of 
the hypersurface $\{x_1=0\}$.
Then it follows by  the monotonicity and  additivity properties of the capacity (see, for example, 
 \cite{FOT}, Section~2.1, or  \cite{BH}, Section~1.8) that  $C_{\overline h}(\{x_1=0\})=0$.
 This establishes the second statement of the proposition.
 
 Finally suppose $\delta\in[0,1/2\rangle$.
 Then one can find a Fourier multiplier $F$ such that $H\geq F$ and, by Proposition~\ref{pcsg3.1},
 one may choose $F$ such that $1+F(p_1,p_2)\geq1+ a\,|p_1|^{2(1-\delta_1)}$ for some $a>0$
 and all $(p_1,p_2)\in \Ri\times \Ri^m$.
 Then if $U\subset \Ri\times\Ri^m$ is an open set with $|U|<\infty$ it follows by the calculation at the end
 of Section~3 of \cite{RSi} that 
 \[
 C_{\overline h}(U)\geq |U|^2(\one_U,(I+F)^{-1}\one_U)^{-1}
 \;.
 \]
 Now set $U_\varepsilon=\langle-\varepsilon,\varepsilon\rangle\times V$ where $V$ is an open subset of $\Ri^m$.
 Then calculating as in the proof of Proposition~4.1 of \cite{RSi} one finds
 \begin{eqnarray*}
  C_{\overline h}(U_\varepsilon)&\geq& 4\,\varepsilon^2\,|V|^2\,\bigg(4\,\varepsilon^2\int_{\Ri^m}dp_2\,(\tilde\one_V(p_2))^2
  \int_\Ri dp_1\,(1+F(p_1,p_2))^{-1}(\sin(\varepsilon p_1)/(\varepsilon p_1))^2\bigg)^{-1}\\[5pt]
  &\geq& |V|\,\bigg(\int_\Ri dp_1\,\left(1+a\,|p_1|^{2(1-\delta_1)}\right)^{-1}\bigg)^{-1}\geq a_{\delta_1}\,|V|
  \end{eqnarray*}
 where $a_{\delta_1}>0$.
 Note that the strict positivity of $a_{\delta_1}$ requires $\delta_1\in[0,1/2\rangle$.
Then,  however, one must have  $C_{\overline h}(\{x_1=0\})>0$.\hfill$\Box$

\bigskip

The moral of the proposition  is that 
in the strongly degenerate case $\delta_1\in[1/2,1\rangle$ the Dirichlet and Neumann boundary conditions coincide.
The separation is a spontaneous effect which is not characterized by a particular choice of boundary conditions.



 \section{A one-dimensional example}\label{Scsg7}
 
In this section we give a further analysis of the one-dimensional example
 discussed in \cite{ERSZ1}, Sections~5 and 6 (see also \cite{ACF} and \cite{MaV} where a similar
 example is analyzed from the point of view of control theory).
 This example is a special case of the Gru\v{s}in operator with $n=1$,  $m=0$ and $\delta_1'=0$.
 Its structure provides a guide to the anticipated structure of the more interesting examples with
  $n=1$ and $m\geq1$.

Let   $h(\varphi)=(\varphi',c_\delta\,\varphi')$ be the closable form on $L_2(\Ri)$ 
 with domain $W^{1,2}(\Ri)$ where $c_\delta$ is given by $c_\delta(x)=(x^2/(1+x^2))^\delta$ with $\delta>0$.
 Let $H$ be  the  positive self-adjoint associated with the 
   closure,  $S$ the semigroup generated by $H$  and $K$ the   kernel of $S$.
 All the qualitative features we subsequently derive extend to the semigroups associated with forms
 $h(\varphi)=(\varphi',c\,\varphi')$ with $c\in L_\infty(\Ri)$ and $c\sim c_\delta$.
 
The  Riemannian distance is now  given by
$d(x\,;y)=|\int^x_ydt\,c_\delta(t)^{-1/2}|$.
If $\delta\in[0,1\rangle$ then $d(\cdot\,;\cdot)$ is a genuine distance but if $\delta\geq 1$
 then the distance between the left and half right lines is infinite.
We concentrate on the  case $\delta\in[0,1\rangle$   and comment on the distinctive features of the case $\delta\geq 1$
at the end of the section.

If  $\delta\in[0,1\rangle$ the Riemannian ball $B(x\,;r)=\{y:d(x\,;y)<r\}$ is an interval
and  the volume is the length of the interval.
It is straightforward to estimate this length from the explicit form of $c_\delta$.
One finds
$|B(x\,;r)|\sim  r $  if  $|x|\geq 1$ or if $r\geq 1$,
$ |B(x\,;r)|\sim r^{1/(1-\delta)}$  if $|x|\leq 1$, $r\leq 1$ and $d(0\,;x)<r$ and 
$ |B(x\,;r)|\sim |x|^\delta\, r$
if $|x|\leq 1$, $r\leq 1$ and $d(0\,;x)\geq r$.
Then $|B|$ satisfies the doubling property (\ref{epre1.17})
with doubling dimension  ${\tilde D} =D=1/(1-\delta)$.
These bounds are all consistent with the general estimates of Proposition \ref{pcsg5.1}.

Moreover, if $B_+(x\,;r)=B(x\,;r)\cap \Ri_+$   then $|B(x\,;r)|/2\leq |B_+(x\,;r)|\leq |B(x\,;r)|$
  for $x\geq0$. 
  Hence $|B_+|$ satisfies similar estimates for $x\geq0$.

\smallskip

Now we consider bounds on the associated semigroup kernel $K_t$.
There are two distinct cases $\delta\in[0,1/2\rangle$ and $\delta\in[1/2,1\rangle$.

\smallskip

\noindent{\bf Case I} $\;\delta\in[0,1/2\rangle$.
In this case the Gaussian upper bounds of Corollary~\ref{ccsg6.1}
 are valid but there are matching lower bounds.

\begin{prop}\label{pcsg4.1}
If $\delta\in[0,1/2\rangle$  there are $b,c>0$ such that 
\begin{equation}
K_t(x\,;y)\geq b\,|B(x\,;t^{1/2})|^{-1}e^{-cd(x;y)^2/t}
\label{ep2}
\end{equation}
for all $x,y\in\Ri$ and $t>0$.
\end{prop}

The deduction of the  off-diagonal lower bounds requires  some additional information
on continuity.

\begin{lemma}\label{lcsg4.4}
If $\delta\in[0,1/2\rangle$ then there is an $a>0$ such that
\[
|\varphi(x)-\varphi(y)|^2\leq a\,{{d(x\,;y)^2}\over{V(x\,;y)}}\,h(\varphi)
\]
for all $x,y\in\Ri$ and all $\varphi\in D(h)$ where $V(x\,;y)=|B(x\,;d(x\,;y))|\vee |B(y\,;d(x\,;y))|$.
\end{lemma}
\proof\
Let $\varphi\in W^{1,2}(\Ri)$.
Then
\begin{eqnarray*}
|\varphi(x)-\varphi(y)|^2&=&\Big|\int^y_xds\,\varphi'(s)\Big|^2\\[5pt]
   &\leq&\Big|\int^y_xds\,c_\delta(s)^{-1}\Big|\,\Big|\int^y_xds\,c_\delta(s)\,\varphi'(s)^2\Big|
     \leq\Big|\int^y_xds\,c_\delta(s)^{-1}\Big|\,h(\varphi)
     \;.
     \end{eqnarray*}
Now the proof of the lemma follows from the upper bound in the next lemma.

\smallskip

Note that at this point it is essential that $\delta\in[0,1/2\rangle$ to ensure that $c_\delta^{-1}$ 
is locally integrable.

\begin{lemma}\label{lcsg4.5}
If $\delta\in[0,1/2\rangle$ then there is an $a>0$ such that
  \[
{{d(x\,;y)^2}\over{|x-y|}}\leq \Big|\int^y_xds\,c_\delta(s)^{-1}\Big|
\leq a\,{{d(x\,;y)^2}\over{V(x\,;y)}}
\]
for all $x,y\in\Ri$.
\end{lemma}
\proof\
The left hand bound follows directly from the Cauchy--Schwarz inequality;
\[
d(x\,;y)^2=\Big|\int^y_xds\,c_\delta(s)^{-1/2}\Big|^2
\leq |x-y|\,\Big|\int^y_xds\,c_\delta(s)^{-1}\Big|
\;.
\]
The right hand bound uses the volume estimates and follows by treating various different cases.

First  set $D(x\,;y)=|\int^y_xds\,c_\delta(s)^{-1}|$.
Then if $|x-y|\geq 1/2$ one has $d(x\,;y)\sim |x-y|$ and by similar reasoning
 $D(x\,;y)\sim |x-y|$.
Moreover, $V(x\,;y)\sim|x-y|$ by  the volume estimate with $r\geq 1$.
Therefore  $D(x\,;y)\sim d(x\,;y)^2/V(x\,;y)$.
If, however, $x,y\geq 1/2$ or $x,y\leq-1/2$ then the estimate follows by similar reasoning
but using the volume estimate with $|x|\geq 1$.
It remains to consider  $x,y$ such that $|x-y|\leq 1$.
But  by symmetry the discussion can be reduced to
 two cases $0\leq y<x\leq 1$ and $-1\leq y<0<x\leq 1$.

 Consider the first case. 
 Then $D(x\,;y)\sim x^{1-2\delta}-y^{1-2\delta}$ and 
$d(x\,;y)\sim x^{1-\delta}-y^{1-\delta}$ by explicit calculation.
 Moreover,  $V(x\,;y)=|B(x\,;d(x\,;y))|\leq a\,(x-y)$ for some $a\geq1$.
 But
\begin{eqnarray*}
(x^{1-\delta}-y^{1-\delta})^2-(x^{1-2\delta}-y^{1-2\delta})(x-y)
=(x^{1/2-\delta}y^{1/2}-x^{1/2}y^{1/2-\delta})^2\geq0
\;.
\end{eqnarray*}
This again establishes the required bound.

Finally consider the second case and suppose that $x\geq |y|$. 
 Then $d(x\,;y)\sim d(x\,;0)\sim x^{1-\delta}$ and $D(x\,;y)\sim D(x\,;0)\sim x^{1-2\delta}$.
 Moreover,
 \[
 V(x\,;y)=|B(x\,;d(x\,;y))|\leq |B(x\,;2d(x\,;0))|\leq a\,|B(x\,;d(x\,;0))|\leq a'\,x
 \]
 where the second estimate uses volume doubling.
 The required bound follows immediately.
 The case $|y|>x$ is similar with the roles of $x$ and $y$ interchanged.
 \hfill$\Box$
 
 \bigskip

The  bound in Lemma~\ref{lcsg4.4} is now an immediate consequence of the upper bound of Lemma~\ref{lcsg4.5}.
Moreover the Gaussian lower bound in Proposition~\ref{pcsg4.1}
 follows directly from the  the Gaussian upper bound and the continuity bound of Lemma~\ref{lcsg4.4}. This last implication is, for example, a direct consequence of Theorem~3.1 in \cite{Cou4}
 applied with $w=2$, $p=2$ and  $\alpha=1$. 
 Note that the doubling dimension
$D=1/(1-\delta)< \alpha \,p=2$ because $\delta\in[0,1/2\rangle$.
Therefore Theorem~3.1 of \cite{Cou4} is indeed applicable.
\hfill$\Box$
 
\bigskip

\noindent{\bf Case 2} $\;\delta\in[1/2,1\rangle$.
It follows from Proposition~6.5 of \cite{ERSZ1} that $S_tL_2(\Ri_\pm)\subseteq L_2(\Ri_\pm)$.
Thus the system separates into two ergodic components and the semigroup kernel has the property $K_t(x\,;y)=0$
for $x<0$ and $y>0$.
One can, however, extend the foregoing analysis to each component.
First we prove that the  kernel of the semigroup restricted to $L_2(\Ri_+)$
is  H\"older continuous.
Note that the generator of the restriction of the semigroup to $L_2(\Ri_+)$ is the operator
associated with the  closure of the form  obtained by restricting $\overline h$  to 
$W^{1,2}(\Ri_+)=\one_{\Ri_+}W^{1,2}(\Ri)$.
Now the continuity  proof is by a variation of the usual   Sobolev inequalities
\[
|\varphi(x)|^2 \leq a\,(\|\varphi'\|_2^2 +\|\varphi\|_2^2)
\]
and
\[
|x-y|^{-2\gamma}|\varphi(x)-\varphi(y)|^2 \leq a_\gamma \,(\|\varphi'\|_2^2 +\|\varphi\|_2^2)
\]
where  $\gamma\in\langle0,1/2\rangle $.
Fix $\sigma>0$ and let  $\psi \in C_c^\infty(\sigma /2,\infty)$ with 
$\psi(x)=1$ for $x>\sigma $.
Then
\begin{eqnarray*}
\|(\psi \varphi)'\|_2^2 +\|\psi \varphi\|_2^2&\leq&
\|\psi' \varphi+\varphi'\psi\|_2^2+\|\psi \varphi \|_2^2 \leq 2\,\|\psi' \varphi\|_2^2+
2\,\|\psi\varphi'\|_2^2 +\|\psi\varphi\|_2^2\\[5pt]
&\leq& 2\,\|\psi\varphi'\|_2^2+(1+4/\sigma )\,\| \varphi\|_2^2\leq 2\,
c_\delta(\sigma /2)^{-1}\|c_\delta^{1/2}\varphi'\|_2^2+(1+4/\sigma )\,\| \varphi\|_2^2\\[5pt]
&\le& a_\sigma \,(h(\varphi)+\|\varphi\|_2^2)\;.
\end{eqnarray*}
Therefore
\[
|\varphi(x)|^2=|(\psi\varphi)(x)|^2 \leq a'_\sigma \,(h(\varphi)+\|\varphi\|_2^2)
\]
and
\[
|x-y|^{-2\gamma}|\varphi(x)-\varphi(y)|^2 \leq a_{\sigma ,\gamma} \,(h(\varphi)+\|\varphi\|_2^2)
\]
for all $x, y\geq \sigma $ and all $\varphi\in W^{1,2}(\Ri_+)$.
These bounds then extend by continuity to the closure of  $h$ and are sufficient to deduce
that the semigroup kernel $K_t$ is uniformly bounded and H\"older continuous on 
$[\sigma ,\infty\rangle\times [\sigma ,\infty\rangle$ for each $\sigma >0$.
In particular $K_t$ is continuous on $\langle 0,\infty\rangle\times \langle 0,\infty\rangle$.
But the bounds depend on $\sigma $ and  do not give good {\it a priori} bounds on 
$\sup_{x,y>0}K_t(x\,;y)=\|S_t\|_{1\to\infty}$.
This can again be accomplished by a Nash inequality argument.

If $h_+$ temporarily denotes the closed form of the generator of the semigroup on $L_2(\Ri_+)$ and 
$E\varphi$ denotes the symmetric extension of $\varphi\in L_2(\Ri_+)$ to $E\varphi\in L_2(\Ri)$ then 
${\overline h}(E\varphi)=2\,h_+(\varphi)$ for all $\varphi\in D(h_+)$.
Moreover, $\|E\varphi\|_2^2=2\,\|\varphi\|_2^2$ and $\|E\varphi\|_1^2=4\,\|\varphi\|_1^2$.
Therefore the Nash inequalities (\ref{epre1.7}) give similar inequalities
\[
\|\varphi\|_2^2\leq r^{-2}h_+(\varphi)
+\pi^{-1}\, V_F(r)\,\|\varphi\|_1^2
\]
for all $\varphi\in D(h_+)$.
Hence the  semigroup restricted to $L_2(\Ri_+)$, or  $L_2(\Ri_-)$, again satisfies bounds
 $\|S_t\|_{1\to\infty}\leq a\,t^{(-1/(2(1-\delta)),-1/2)}$ for all $t>0$.
In particular, the kernel $K_t$ is uniformly bounded on $\Ri_+\times \Ri_+$,
or  $\Ri_-\times \Ri_-$,
for each $t>0$.

One can then prove the analogues of Corollaries~\ref{ccsg6.1}
and ~\ref{ccsg6.2} on the half-lines.
\begin{prop}\label{pcsg4.2}
If $\delta\in[1/2,1\rangle$ and  $B_\pm(x\,;r)=B(x\,;r)\cap \Ri_\pm$ then
for each $\varepsilon\in\langle0,1]$ there is an $a>0$ such that 
\begin{equation}
K_t(x\,;y)\leq a\,|B_\pm(x\,;t^{1/2})|^{-1}e^{-d(x;y)^2/(4t(1+\varepsilon))}
\label{ep3}
\end{equation}
for all $x,y\in\Ri_\pm$ and $t>0$.
Moreover, there is  $b>0$ such that 
\begin{equation}
K_t(x\,;x)\geq b\,|B_\pm(x\,;t^{1/2})|^{-1}\label{ep4}
\end{equation}
for all $x\in\Ri_\pm$ and $t>0$.
\end{prop}

In the one-dimensional case the statement of Proposition~\ref{psep}
 can also be described in terms of the vector fields defining $h$ and $h_{\rm Dir}$.

 Let $X=c_\delta^{1/2}\,d$, with $d=d/dx$,  denote the $C^\delta$-vector field acting on $L_2(\Ri)$  with domain $D(X)=C_c^\infty(\Ri)$  and $X_0$ the restriction of $X$ to 
 $D(X_0)=C_c^\infty(\Ri\backslash\{0\})$.
 Then $h(\varphi)=\|X\varphi\|_2^2$ and $h_{\rm Dir}(\varphi)=\|X_0\varphi\|_2^2$.
 It follows straightforwardly that the corresponding positive self-adjoint operators are given by
 $H=X^*\overline X$ and $H_{\rm Dir}=X_0^*\overline X_0$.
 But Proposition~\ref{psep} can now be restated as follows.
 
 \begin{cor}\label{csep}
    If $\delta\in[0,1/2\rangle$ then  $\overline X_0\subset \overline X$ but if  $\delta\in[1/2,1\rangle$ then  $\overline X_0= \overline X$.
    \end{cor}

Although we have restricted attention to the case $\delta\in[0,1\rangle$ the form $h$ is densely defined
and closable for all $\delta\geq 0$.
But the semigroup and its kernel have very different properties if  $\delta\geq 1$.
The properties of the distance
$d(x\,;y)=|\int^x_ydt\,c_\delta(t)^{-1/2}|$ are also quite different.
The distance is finite on the open half line $\langle 0,\infty\rangle$ but since the integral diverges at zero
the distance to the origin is at infinity.
Now consider points $x,y\in\langle0,1]$. Then $d(x\,;y)\sim |\ln x/y\,|$ if $\delta=1$ and $d(x\,;y)\sim|x^{1-\delta}-y^{1-\delta}|$
if $\delta>1$.
Consider the case $\delta=1$ with the equivalent distance $\tilde d(x\,;y)=|\ln x/y\,|$. 
Let $x_n\in\langle0,1\rangle$ be an arbitrary sequence which converges downward to zero.
One may assume $x_1\leq e^{-1}$.
Then $\tilde d(x_n\,;x_ne)=1=\tilde d(x_n\,;x_ne^{-1})$ and $|B(x_n\,;1)|=2\,x_n\, \sinh 1$.
Therefore $|B(x_n\,;1)|\to 0$ as $n\to\infty$.

Alternatively if $x_n=e^{-n}$ with $n\geq 1$ then $B(x_n\,;n/2)=\langle e^{-3n/2},e^{-n/2}\rangle$ and 
$B(x_n\,;n)=\langle e^{-2n},1\rangle$. 
Therefore $|B(x_n\,;n/2)|=e^{-n/2}(1-e^{-n})\to0$ and  
$|B(x_n\,;n)|=1-e^{-2n}\to1$ as $n\to\infty$ so the volume cannot satisfy the volume
doubling property.

These divergences allow one to argue that the semigroup kernel is not bounded near the origin.

\smallskip


%


 \section{Applications}\label{Scsg8}
 
 In the foregoing we  established that Gru\v{s}in operators have many important properties
 in common with strongly elliptic operators; the wave equation has a finite propagation speed,
 the heat kernel satisfies Gaussian upper bounds and the heat semigroup conserves probability.
 Then one can readily adapt arguments developed for strongly elliptic operators to obtain
 further detailed information about the Gru\v{s}in operators, e.g.\ information on boundedness of Riesz 
 transforms, spectral multipliers and Bochner--Riesz summability, holomorphic functional calculus, Poincar\'e inequalities and maximal   regularity.
 We conclude by describing briefly some of these applications to Gru\v{s}in operators which are 
 a straightforward consequence of general theory and which require no further detailed arguments.
 It should, however, be emphasized that there are significant differences between the degenerate and the non-degenerate  theories  related to continuity and positivity properties.
 In particular the one-dimensional example in Section~\ref{Scsg7} demonstrates that 
 the heat kernel  is not necessarily continuous nor strictly positive.
 Therefore there are limitations to possible extensions of the results of classical analysis to the degenerate case.

\subsection{Boundedness of Riesz transforms}

First we consider boundedness of the Riesz  transforms associated with a general Gru\v{s}in operator $H$ on $L_p(\Ri^{n+m})$ for $p\in\langle1,2] $. 
The result can be stated  in terms of the {\bf carr\'e du champ} associated with $H$
(see, for example, Section~I.4 of \cite{BH}). 
Formally  the carr\'e du champ is given by 
$\Gamma_\psi=\psi(H\psi)-2^{-1}H\psi^2$.
Note that if $\psi \in C_c^\infty(\Ri^{n+m})$ then
\[
 \Gamma_\psi(x)=\sum_{i,j=1}^{n+m}c_{ij}(x)(\partial_i \psi)(x)( \partial_j\psi)(x)
\]
and
$\| \Gamma_\psi\|_1=\|H^{1/2}\psi\|_2$.
Now one can formulate the
result concerning boundedness of  the   Riesz
transform in an analogous manner.
\begin{thm}\label{tcsgc.1}
If $H$ is   Gru\v{s}in  operator then
$$
\|\Gamma_{H^{-1/2}\psi}\|_{p/2}\leq \|\psi\|_{p}
$$
for all $\psi \in L^p(\Ri^{n+m})$ and all
$p\in\langle1,2]$.
In addition the map $\psi \to
\Gamma_{H^{-1/2}\psi}^{1/2}$ is  weak type~$(1,1)$, i.e.\
\[
|\{ x \in X : |  \Gamma_{H^{-1/2}\psi} (x)|^{1/2} >
\lambda \}
\leq a \, {\| \psi \|_{1}}/{\lambda}
\]
for all $ \lambda \in \Ri_+$ and all  $ \psi \in L_1(\Ri^{n+m})$.
\end{thm}
\proof\
The proof of Theorem~\ref{tcsgc.1} is a straightforward modification of
the proof of Theorem~5 of \cite{Sik3}.
The assumptions of Theorem~5 of \cite{Sik3} hold in virtue of the property of
finite speed of propagation proved
in Proposition~4.1 and the kernel bounds (\ref{ecsg6.2}) of Theorem~\ref {tcsg6.2} .
\hfill$\Box$

\subsection{Spectral multipliers}
Each  Gru\v{s}in operator  $H$ is  positive definite and
self-adjoint.
Therefore $H$ admits a
spectral  resolution  $E_H(\lambda)$  and for  any  bounded  Borel
function $F\colon [0, \infty) \to \Ci$ one can define the operator $F(H)$
by 
   \begin{equation}\label{equw}
   F(H)=\int_0^{\infty}d E_H(\lambda)\,F(\lambda) \;.
   \end{equation}
It then follows that $F(H)$ is bounded on
$L^2(\Ri^{n+m})$.
Spectral  multiplier  theorems  investigate  sufficient conditions on
function  $F$ which  ensure  that the   operator  $F(H)$ extends  to a
bounded operator on  $L_q$ for some $q\in[1,\infty]$.
\begin{thm}\label{tcsgc.2}
If $H$ is a  Gru\v{s}in  operator,
$s > (D\vee D')/2$ and $F \colon [0,\infty) \to \Ci$
 is  a  bounded Borel function  such that
     \begin{equation}\label{Ale}
     \sup_{t>0}\| \eta \,\delta_t F \|_{W^{s,\infty}} < \infty,
     \end{equation}
   where $ \delta_t F(\lambda)=F(t\lambda)$ and
   $\| F \|_{W^{s,p}}=\|(I-d^2/d x^2)^{s/2}F\|_{L_p}$. 
   Then $F(H)$ is
 weak type $(1,1)$ and bounded on~$L_q$ for all  $q\in\langle1,\infty\rangle$.
\end{thm}
\proof\
The proof of Theorem~\ref{tcsgc.2} is a direct consequence of the Gaussian bounds  (\ref{ecsg6.21}) on the heat kernel corresponding to $H$
given by Corollary~\ref{ccsg6.1} and Theorem~3.1 of \cite{DOS}.
(See also Theorem~3.5 of \cite{CS}.)
\hfill$\Box$

\bigskip

The theory of spectral multipliers is  related to  and motivated
by  the study of convergence  of the Riesz   means or  convergence   of other
eigenfunction expansions of self-adjoint  operators.
To  define  the
Riesz means of the operator $H$ we set
   \begin{equation}\label{vab1}
   \sigma^{s}_R(\lambda)=
      \left\{
       \begin{array}{cl}
       (1-\lambda/R)^{s}  &\mbox{for}\;\; \lambda \le R \\
       0  &\mbox{for}\;\; \lambda > R. \\
       \end{array}
      \right.
   \end{equation}
We then define  $\sigma^{s}_R(H)$  by spectral theory.
The operator $\sigma^{s}_R(H)$ is the Riesz or the Bochner-Riesz mean of
order $s$. 
The basic question in the theory of the Riesz means is  to
establish  the critical exponent for continuity and convergence of
the Riesz means.
More precisely one wishes to ascertain  the optimal range of
$s$ for which the Riesz means $\sigma^{s}_R(H)$ are uniformly
bounded on $L_1(\Ri^{n+m})$. 
A result of this type is given by the
following.
\begin{thm}\label{tcsgc.3}
If $H$ is a  Gru\v{s}in  operator and     $s > (D\vee D')/2$
then
\[
   \sup_{R>0}\|\sigma_R^{s}(H)\|_{q \to q}
     \leq a < \infty
\]
 for all $q\in [1,\infty]$.
   Hence  
\[
     \lim_{R \to \infty}
     \|\sigma_R^{s}(H)\varphi-\varphi\|_{q \to q}=0
\]
for all $q\in [1,\infty]$ and all $\psi\in L_q(\Ri^{n+m})$.
\end{thm}
\proof\
The proof of Theorem \ref{tcsgc.3} is again a  direct consequence of
the Gaussian bounds (\ref{ecsg6.21}) on the heat kernel corresponding to $H$
and
Corollary~6.3 of \cite{DOS}.
\hfill$\Box$

\bigskip

Next we consider the implication of  the Gaussian bounds on the heat kernel
for the holomorphic
function calculus of the Gru\v{s}in operators.
First we briefly recall the notion of holomorphic
function calculus.
For each $\theta >  0$ set  $\Sigma(\theta)=\{z\in C\backslash\{0\}\colon
|\mbox{arg}\,z|<\theta$\}. Let $F$ be a
bounded   holomorphic  function  on     $\Sigma(\theta).$
By   $\Vert F \Vert_{\theta,\infty}$ we denote the supremum of $F$ on
$\Sigma(\theta)$.
The general problem of interest is to  find sharp bounds, in terms of
$\theta$, of the norm of $F(H)$ as an  operator  acting on
$L_p(\Ri^{n+m})$.
  It is known (see  \cite{CDMY}, Theorem~4.10) that
these bounds on  the   holomorphic
functional calculus when $\theta$ tends to~$0$ are related to
spectral multiplier theorems for $H$. 
The following theorem describing
  holomorphic function calculus for  Gru\v{s}in  type operators
follows from (\ref{ecsg6.21}) and Corollary~\ref{ccsg6.1}.
\begin{thm}\label{tcsgc.4}
If  $H$ is a  Gru\v{s}in  operator and 
$s> (D\vee D')|1/p-1/2|$
then
\[
\|  F(H) \|_{p \to p} \le {a}\,{\theta^{-s}}
\| F  \|_{\theta, \infty}
\]
for all $\theta > 0$.
\end{thm}
\proof\
Theorem \ref{tcsgc.4} follows from
the Gaussian bounds (\ref{ecsg6.21}) of
Corollary~\ref{ccsg6.1} and
Proposition~8.1 of \cite{DOS}.
\hfill$\Box$

\subsection{Concluding remarks and comments}

The above statements on the boundedness of the Riesz transforms
and the spectral multipliers for Gru\v{s}in  operators  are not always optimal. 
Using the basic estimates of 
Corollary~\ref{ccsg6.1} and Proposition~4.1 one can analyze the
boundedness
of the Riesz transforms for $p>2$.

In the multiplier result  discussed above, Theorem~\ref{tcsgc.2}, the
critical
exponent required for the  order of differentiability of the function
$F$ is equal to  half of the homogeneous dimension $D\vee D'$.
This  is a  quite typical situation and  for the standard Laplace operator
this exponent is optimal. 
We expect,  however,  that in many cases  it is
possible to obtain multiplier results for Gru\v{s}in  operators with
critical exponent essentially smaller then the half of the homogeneous
dimension $D\vee D'$.

It is also possible to obtain a  version of the Poincare inequality and Nash type
results
similar to those discussed in Section~\ref{Scsg7}. 
But  results of this nature require
substantial new proofs which we hope to describe elsewhere
We conclude by  stressing that
Corollary~\ref{ccsg6.1} and Proposition~\ref{pscsg4.0}  provide a sound basis for
further analysis of Gru\v{s}in
type operators.

\subsection*{Acknowledgement}
This work  was  supported by an Australian Research Council (ARC) Discovery Grant DP 0451016.
It grew out of an earlier collaboration with Tom ter Elst to whom the authors are indebted for many
helpful discussions about degenerate operators.
The work was completed whilst the first author was a guest of Prof.\ Ola Bratteli at the University of Oslo.


\newpage


\begin{thebibliography}{CDMY96}

\bibitem[ABCF06]{ACF}
{\sc Alabau-Boussouira, F., Cannarsa, P., {\rm and} Fragnelli, G.}, Carleman
  estimates for degenerate parabolic operators with applications to null
  controllability.
\newblock {\em J. Evol.\ Equ.} {\bf 6}, No.\ 2 (2006),  161--204.

\bibitem[AH05]{AH}
{\sc Ariyoshi, T., {\rm and} Hino, M.}, Small-time asymptotic estimates in
  local Dirichlet spaces.
\newblock {\em Elec.\ J. Prob.} {\bf 10} (2005),  1236--1259.

\bibitem[BH91]{BH}
{\sc Bouleau, N., {\rm and} Hirsch, F.}, {\em Dirichlet forms and analysis on
  Wiener space}, vol.\ 14 of de Gruyter Studies in Mathematics.
\newblock Walter de Gruyter \& Co., Berlin, 1991.

\bibitem[BR97]{BR2}
{\sc Bratteli, O., {\rm and} Robinson, D.~W.}, {\em Operator algebras and
  quantum statistical mechanics}, vol.\ 2.
\newblock Second edition. Springer-Verlag, New York etc., 1997.

\bibitem[Bra02]{Bra}
{\sc Braides, A.}, {\em $\Gamma$-convergence for beginners}, vol.\ 22 of Oxford
  Lecture Series in Mathematics and its Applications.
\newblock Oxford University Press, Oxford, 2002.

\bibitem[CDMY96]{CDMY}
{\sc Cowling, M., Doust, I., McIntosh, A., {\rm and} Yagi, A.}, Banach space
  operators with a bounded $H^\infty$ functional calculus.
\newblock {\em J. Austr.\ Math.\ Soc.\ {\rm (}Series A{\rm )}} {\bf 60} (1996),
   51--89.

\bibitem[CKN84]{CKN}
{\sc Caffarelli, L., Kohn, R., {\rm and} Nirenberg, L.}, First order
  interpolation inequalities with weights.
\newblock {\em Compositio Math.} {\bf 53} (1984),  259--275.

\bibitem[CKS87]{CKS}
{\sc Carlen, E.~A., Kusuoka, S., {\rm and} Stroock, D.~W.}, Upper bounds for
  symmetric Markov transition functions.
\newblock {\em Ann.\ Inst.\ Henri Poincar\'e} {\bf 23} (1987),  245--287.

\bibitem[Cou03]{Cou4}
{\sc Coulhon, T.}, Off-diagonal heat kernel lower bounds without Poincar\'e.
\newblock {\em J. London Math.\ Soc.} {\bf 68} (2003),  795--816.

\bibitem[CSC95]{CS}
{\sc Coulhon, T., {\rm and} Saloff-Coste, L.}, Vari\'et\'es riemanniennes
  isom\'etriques \`a l'infini.
\newblock {\em Rev.\ Mat.\ Iberoamericana} {\bf 11} (1995),  687--726.

\bibitem[{Dal}93]{DalM}
{\sc {Dal Maso}, G.}, {\em An introduction to {$\Gamma$}-convergence}, vol.\ 8
  of Progress in Nonlinear Differential Equations and their Applications.
\newblock Birkh{\"a}user Boston Inc., Boston, MA, 1993.

\bibitem[Dav92]{Dav12}
{\sc Davies, E.~B.}, Heat kernel bounds, conservation of probability and the
  Feller property.
\newblock {\em J. Anal.\ Math.} {\bf 58} (1992),  99--119.
\newblock Festschrift on the occasion of the 70th birthday of Shmuel Agmon.

\bibitem[Dav99]{Dav13}
\leavevmode\vrule height 2pt depth -1.6pt width 23pt, A review of Hardy
  inequalities.
\newblock In {\em The Maz'ya anniversary collection, Vol.\ 2 (Rostock, 1998)},
  vol.\ 110 of Oper.\ Theory Adv.\ Appl.,  55--67. Birkh{\"a}user, Basel, 1999.

\bibitem[DOS02]{DOS}
{\sc Duong, X.~T., Ouhabaz, E.-M., {\rm and} Sikora, A.}, Plancherel-type
  estimates and sharp spectral multipliers.
\newblock {\em J. Funct.\ Anal.} {\bf 196} (2002),  443--485.

\bibitem[ERS07]{ERS4}
{\sc Elst, A. F.~M. ter, Robinson, D.~W., {\rm and} Sikora, A.}, Small time
  asymptotics of diffusion processes.
\newblock {\em J. Evol.\ Equ.} {\bf 7} (2007),  79--112.

\bibitem[ERSZ07]{ERSZ1}
{\sc Elst, A. F.~M. ter, Robinson, D.~W., Sikora, A., {\rm and} Zhu, Y.},
  Second-order operators with degenerate coefficients.
\newblock {\em Proc.\ London Math.\ Soc.} (2007).
\newblock To appear, doi: 10.1112/plms/pdl017.

\bibitem[ET76]{ET}
{\sc Ekeland, I., {\rm and} Temam, R.}, {\em Convex analysis and variational
  problems}.
\newblock North-Holland Publishing Co., Amsterdam, 1976.

\bibitem[FGW94]{FGW2}
{\sc Franchi, B., Guti{\'e}rrez, C.~E., {\rm and} Wheeden, R.~L.}, Weighted
  Sobolev--Poincar\'e inequalities for Grushin type operators.
\newblock {\em Comm.\ Part.\ Diff.\ Eq.} {\bf 19} (1994),  523--604.

\bibitem[FKS82]{FKS}
{\sc Fabes, E.~B., Kenig, C.~E., {\rm and} Serapioni, R.~P.}, The local
  regularity of solutions of degenerate elliptic equations.
\newblock {\em Comm.\ Part.\ Diff.\ Eq.} {\bf 7} (1982),  77--116.

\bibitem[FL83]{FL}
{\sc Franchi, B., {\rm and} Lanconelli, E.}, H{\"o}lder regularity theorem for
  a class of linear nonuniformly elliptic operators with measurable
  coefficients.
\newblock {\em Ann.\ Scuola Norm.\ Sup.\ Pisa Cl.\ Sci.} {\bf 10} (1983),
  523--541.

\bibitem[FL84]{FL2}
\leavevmode\vrule height 2pt depth -1.6pt width 23pt, An embedding theorem for
  Sobolev spaces related to nonsmooth vector fields and Harnack inequality.
\newblock {\em Comm.\ Part.\ Diff.\ Eq.} {\bf 9} (1984),  1237--1264.

\bibitem[FOT94]{FOT}
{\sc Fukushima, M., Oshima, Y., {\rm and} Takeda, M.}, {\em Dirichlet forms and
  symmetric Markov processes}, vol.\ 19 of de Gruyter Studies in Mathematics.
\newblock Walter de Gruyter \& Co., Berlin, 1994.

\bibitem[Fra91]{Fra}
{\sc Franchi, B.}, Weighted Sobolev-Poincar\'e inequalities and pointwise
  estimates for a class of degenerate elliptic equations.
\newblock {\em Trans.\ Amer.\ Math.\ Soc.} {\bf 327} (1991),  125--158.

\bibitem[FS87]{FrS}
{\sc Franchi, B., {\rm and} Serapioni, R.}, Pontwise estimates for a class of
  strongly degenerate elliptic operators: a geometrical approach.
\newblock {\em Ann.\ Scuola Norm.\ Sup.\ Pisa Cl.\ Sci.} {\bf 14} (1987),
  527--568.

\bibitem[Gaf59]{Gaf}
{\sc Gaffney, M.~P.}, The conservation property of the heat equation on
  Riemannian manifolds.
\newblock {\em Comm.\ Pure Appl.\ Math.} {\bf 12} (1959),  1--11.

\bibitem[Gri99]{Gri3}
{\sc Grigor'yan, A.}, Estimates of heat kernels on Riemannian manifolds.
\newblock In {\em Spectral theory and geometry $($Edinburgh, 1998$)$}, vol.\
  273 of London Math.\ Soc.\ Lecture Note Ser.,  140--225. Cambridge Univ.\
  Press, Cambridge, 1999.

\bibitem[Gru70]{Gru}
{\sc Gru\v{s}in, V.~V.}, A certain class of hypoelliptic operators.
\newblock {\em Mat.\ Sb.\ (N.S.)} {\bf 83 $($125$)$} (1970).

\bibitem[GT83]{GT}
{\sc Gilbarg, D., {\rm and} Trudinger, N.~S.}, {\em Elliptic partial
  differential equations of second order}.
\newblock Second edition, Grundlehren der mathematischen Wissenschaften 224.
  Springer-Verlag, Berlin etc., 1983.

\bibitem[Haj96]{Haj}
{\sc Hajlasz, P.}, Geometric approach to Sobolev spaces and badly degenerated
  elliptic equations.
\newblock In {\em Nonlinear analysis and applications (Warsaw, 1994)}, vol.\ 7
  of Gakuto Internat.\ Ser.\ Math.\ Sci.\ Appl.,  141--168. Gakk\=otosho,
  Tokyo, 1996.

\bibitem[HR03]{HiR}
{\sc Hino, M., {\rm and} Ram\'irez, J.~A.}, Small-time Gaussian behavior of
  symmetric diffusion semigroups.
\newblock {\em Ann.\ Prob.} {\bf 31} (2003),  254--1295.

\bibitem[Jos98]{Jos}
{\sc Jost, J.}, Nonlinear Dirichlet forms.
\newblock In {\em New directions in Dirichlet forms}, vol.\ 8 of AMS/IP Stud.\
  Adv.\ Math.,  1--47. Amer.\ Math.\ Soc., Providence, RI, 1998.

\bibitem[JSC87]{JSC1}
{\sc Jerison, D., {\rm and} S{\'a}nchez-Calle, A.}, Subelliptic, second order
  differential operators.
\newblock In {\sc Berenstein, C.~A.}, ed., {\em Complex analysis III}, Lecture
  Notes in Mathematics 1277. Springer-Verlag, Berlin etc., 1987,  46--77.

\bibitem[Mos94]{Mosco}
{\sc Mosco, U.}, Composite media and asymptotic Dirichlet forms.
\newblock {\em J. Funct.\ Anal.} {\bf 123} (1994),  368--421.

\bibitem[MR92]{MR}
{\sc Ma, Z.~M., {\rm and} R{\"o}ckner, M.}, {\em Introduction to the theory of
  (non symmetric) Dirichlet Forms}.
\newblock Universitext. Springer-Verlag, Berlin etc., 1992.

\bibitem[MV06]{MaV}
{\sc Martinez, P., {\rm and} Vancostenoble, J.}, Carleman estimates for
  one-dimensional degenerate heat equations.
\newblock {\em J. Evol.\ Equ.} {\bf 6} (2006),  325--362.

\bibitem[Nas58]{Nash}
{\sc Nash, J.}, Continuity of solutions of parabolic and elliptic equations.
\newblock {\em Amer.\ J. Math.} {\bf 80} (1958),  931--954.

\bibitem[Rob91]{Robm}
{\sc Robinson, D.~W.}, {\em Elliptic operators and Lie groups}.
\newblock Oxford Mathematical Monographs. Oxford University Press, Oxford etc.,
  1991.

\bibitem[RS05]{RSi}
{\sc Robinson, D.~W., {\rm and} Sikora, A.}
\newblock Degenerate elliptic operators: capacity, flux and separation, 2005.
\newblock arXiv:math.AP/0601351.

\bibitem[SC02]{Sal5}
{\sc Saloff-Coste, L.}, {\em Aspects of Sobolev-type inequalities}.
\newblock London Math.\ Soc.\ Lect.\ Note Series 289. Cambridge University
  Press, Cambridge, 2002.

\bibitem[Sik96]{Sik}
{\sc Sikora, A.}, Sharp pointwise estimates on heat kernels.
\newblock {\em Quart.\ J. Math.\ Oxford} {\bf 47} (1996),  371--382.

\bibitem[Sik04]{Sik3}
\leavevmode\vrule height 2pt depth -1.6pt width 23pt, Riesz transform, Gaussian
  bounds and the method of wave equation.
\newblock {\em Math.\ Z.} {\bf 247} (2004),  643--662.

\bibitem[Str67]{Stri}
{\sc Strichartz, R.~S.}, Multipliers on fractional Sobolev spaces.
\newblock {\em J. Math.\ Mech.} {\bf 16} (1967),  1031--1060.

\bibitem[SW06]{SawW}
{\sc Sawyer, E.~T., {\rm and} Wheeden, R.~L.}, H{\"o}lder continuity of weak
  solutions to subelliptic equations with rough coefficients.
\newblock {\em Mem.\ Amer.\ Math.\ Soc.} {\bf 180}, No.\ 847 (2006),  x+157.

\bibitem[Tru73]{Tru2}
{\sc Trudinger, N.~S.}, Linear elliptic operators with measurable coefficients.
\newblock {\em Ann.\ Scuola Norm.\ Sup.\ Pisa} {\bf 27} (1973),  265--308.

\end{thebibliography}

\end{document}